\documentclass[journal]{new-aiaa}

\title{String stability of energy-saving aircraft formations}
\author{James R. Riehl \footnote{Research staff, Institute of Mechanics, Materials and Civil Engineering}}
\author{Esteban A. L. Hufstedler\footnote{Postdoctoral researcher, Institute of Mechanics, Materials and Civil Engineering}}
\author{Philippe Chatelain\footnote{Professor, Institute of Mechanics, Materials and Civil Engineering}}
\author{Julien M. Hendrickx\footnote{Professor, Institute of Information and Communication Technologies, Electronics and Applied Mathematics}}
\affil{UCLouvain, 1348 Louvain-la-Neuve, Belgium. This work was 
supported by the “RevealFlight” Concerted Research Action (ARC) of
the Federation Wallonie-Bruxelles}

\usepackage{textcomp}
\usepackage{graphicx}
\usepackage{amsmath}
\usepackage{siunitx}
\usepackage{longtable,tabularx}
\usepackage[textsize=small,disable]{todonotes}
\usepackage[margin=1in]{geometry}
\usepackage{graphicx}
\usepackage{amsmath}
\usepackage[title]{appendix}
\usepackage{tikz}
\usetikzlibrary{shapes,arrows}
\usetikzlibrary{calc,positioning}
\usepackage{siunitx}
\usepackage{color}

\usepackage{caption}
\usepackage{subcaption}
\usepackage{lscape}
\usepackage{ulem}

\newcommand{\e}{\mathbf{e}}
\renewcommand{\v}{\mathbf{v}}

\newcommand{\x}{\mathbf{x}}

\newcommand{\p}{\mathbf{p}}
\renewcommand{\a}{\boldsymbol{\alpha}}

\renewcommand{\u}{\mathbf{u}}

\newcommand{\w}{\mathbf{w}}

\newcommand{\control}{\mathbf{u}}

\newcommand{\pos}{\mathbf{p}}

\newcommand{\sepRef}{\boldsymbol{\delta}_\text{ref}}
\newcommand{\sepRefElem}{\delta}
\newcommand{\err}{\mathbf{e}}

\newcommand{\vel}{\mathbf{v}}
\newcommand{\velElem}{v}

\newcommand{\headway}{h}

\newcommand{\revjr}{}
\newcommand{\revjh}{} %[1]{\textcolor{blue}{#1}}
\newcommand{\rev}{} % \textcolor{blue}}
\newcommand{\revTwo}{} % {\textcolor{green}}
\newcommand{\revThree}{} %{\textcolor{red}} % {\textcolor{green}}
\newcommand{\done}[1]{\todo[disable]{#1}\addcontentsline{tdo}{todo}{\sout{#1}}}

\tikzstyle{int}=[draw, fill=blue!20, minimum size=2em]
\tikzstyle{init} = [pin edge={to-,thin,black}]

\tikzset{
    block/.style = {draw, rectangle,
        minimum height=1cm,
        minimum width=2cm},
    input/.style = {coordinate,node distance=1cm},
    output/.style = {coordinate,node distance=4cm},
    arrow/.style={draw, -latex,node distance=2cm},
    pinstyle/.style = {pin edge={latex-, black,node distance=2cm}},
    sum/.style = {draw, circle, node distance=1cm},
}

\tikzstyle{block} = [draw, fill=blue!20, rectangle, 
    minimum height=3em, minimum width=6em]
\tikzstyle{sum} = [draw, fill=blue!20, circle, node distance=1cm]
\tikzstyle{input} = [coordinate]
\tikzstyle{output} = [coordinate]
\tikzstyle{feedback} = [coordinate]
\tikzstyle{pinstyle} = [pin edge={to-,thin,black}]

\def\code#1{\texttt{#1}}

\begin{document}
\newpage

\maketitle

\begin{abstract}
    Groups of aircraft have the potential to save significant amounts of energy by flying in formations; all but the leading aircraft can benefit from the upwash of the wakes of preceding aircraft.
    A potential obstacle as the number of aircraft in such a formation increases is that disturbances at one aircraft, for example caused by turbulence or wake meandering, can propagate and grow as each following aircraft tries to track the optimal energy-saving position relative to the one in front.
    This phenomenon, known as \textit{string instability}, has not yet been adequately examined \rev{in the context of aircraft formations}.
    We discuss some trade-offs involved in designing string stable controllers whose objective is to minimize \revjh{energy, and present} a control design method to achieve both string stability and energy efficiency of an aircraft formation.
    In simulations of a 10-aircraft linear formation in the presence of 2\% turbulence intensity, \revjr{our controller achieves string stability while reduced energy consumption by an average of 13\% with respect to solo flight.} 
\end{abstract}

\section{Introduction}

\begin{figure}[ht]
    \centering
    \includegraphics[width=0.5\linewidth]{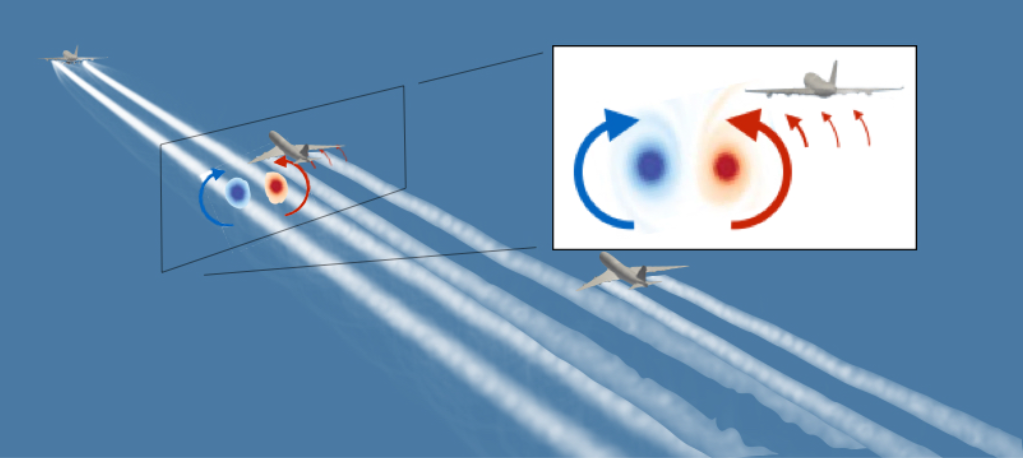}
    \caption{Airplanes in flight produce pairs of trailing vortices that rotate in the directions indicated by the red and blue arrows. An airplane flying behind another that is positioned in the upward flowing region of the vortices can save energy as a result of this additional lift.}
    \label{fig:wakeVelocity}
\end{figure}

There is evidence that migratory bird flocks are able to reduce their energy use by around 12\% by flying in formations that exploit the aerodynamic benefits of wakes \cite{voelkl2015matching}.
Groups of fixed-wing aircraft have the potential to save energy by the same principle, and this has already been demonstrated on pairs of aircraft in wind tunnels \cite{blake2004comparison} as well as in free flight \cite{hansen2002induced,okolo2014effect}.
Briefly, this is possible because winged agents in flight generate pairs of vortices containing regions of upward airflow (upwash) and downward airflow (downwash), as shown in Fig. \ref{fig:wakeVelocity}. 
While simplified wake models can predict the approximate locations of these regions, complex aerodynamic effects such as wake meandering, turbulence, and wake-wake interactions render the task of tracking these locations quite challenging in practice. 
But even supposing the \revjr{location of the wakes can be measured precisely}, there remains the problem that cascaded formations such as these are subject to disturbance \rev{amplification from agent to agent, which can grow unbounded as the number of agents increases.
\rev{This phenomenon, known as \textit{string instability} arises frequently in the literature on vehicle platoons traveling on roads \cite{swaroop1994comparison,ploeg2014lp,farnam2018strong}.}
For example, a small deceleration by one vehicle in a long sequence of self-driving vehicles on a highway might grow from one car to the next and eventually result in a collision.
Although not yet thoroughly analyzed in the context of aircraft formations, the following simple example shows that it can also play a role here.}

\begin{figure}[ht]
    \centering 
    \includegraphics[width=0.5\linewidth]{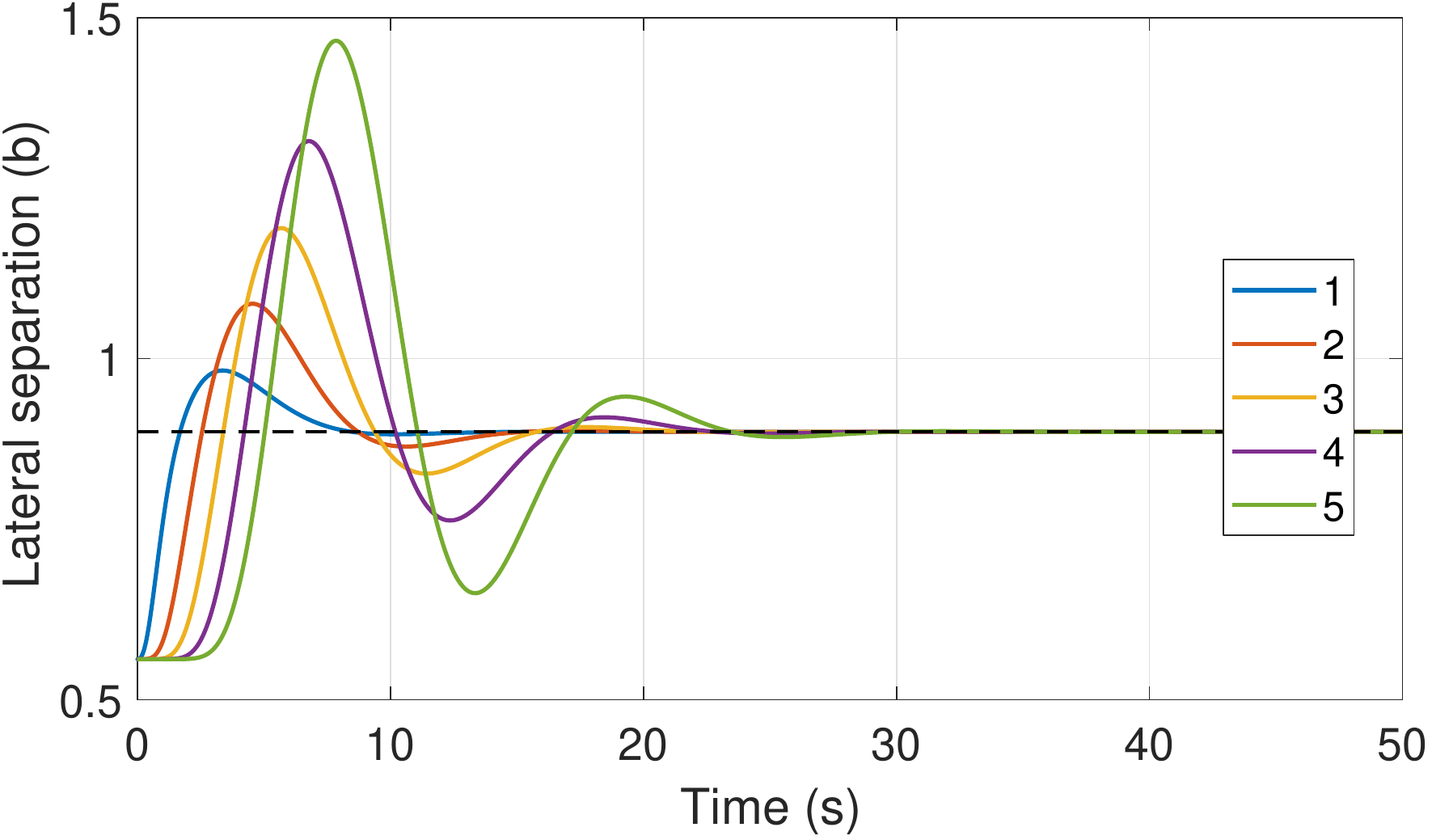}
    \caption{A small disturbance to the lateral position of the formation leader is amplified in the response of each of the four successive following aircraft. \rev{Lateral separation is given in units of wingspan (b). The energy-optimal separation is indicated by the dashed line, and the airplanes start the simulation with a smaller separation. The aircraft model used in the simulation is discussed further in Section \ref{model}.}}
    \label{fig:prev}
\end{figure}

\rev{Consider a linear formation in which the objective of each following aircraft is to maintain a prescribed offset vector with the preceding aircraft. Suppose we attempt this by linearizing around some desired steady-state conditions and applying a classical LQR controller \cite{stevens1992aircraft}.} 
Fig. \ref{fig:prev} shows that when five airplanes are connected in this manner, small disturbances can propagate and grow along the formation resulting in much larger perturbations towards the end.
\rev{Although the simulated aircraft are able to recover in this example, it is easy to see that a large enough perturbation or long enough formation would result in the controller commanding dangerous maneuvers.} \revjh{Moreover, this phenomenon suggests that a small recurrent disturbance at the leader would lead to large permanent effects down the formation.} As we will demonstrate later on, it is possible to design controllers with better disturbance attenuation than shown here, but this requires carefully taking such propagation phenomena into account during the design phase.

\rev{In the rare instances when string stability has been considered in the literature on aircraft formations, it has presented a challenge for the control design.
For example, in \cite{allen2002string}, a PID leader-follower controller for two F/A-18 aircraft was extended to longer formations. While the formation was steady enough to comply with ride quality requirements for up to seven aircraft, it was shown to exhibit string unstable behavior.
The authors suggested that string unstable controllers could be tolerated for limited formation sizes.} \revTwo{However, energy consumption was not a concern in their context.}
\revjr{In \cite{swieringa2015string,weitz2018comparing}, string stability and steady-state error were considered in the context of interval management, where air traffic controllers provide speed guidance to regulate the spacing between aircraft on landing approach.
This demonstrates that even in a centralized control context that does not incorporate aircraft dynamic models or wake effects, string stability can be a relevant factor.
}
\rev{In this paper, we investigate the extent to which string stability is a problem in energy-saving aircraft formations and potential remedies for such problems in the control design.}

\rev{First, we show that in contrast to automobile platoons \revjh{for which fundamental impossibility results have been proved}, string stability of an aircraft formation based on a classical dynamic model with fixed separation distance can be achieved by standard \revjr{state feedback control, e.g. a linear quadratic regulator (LQR).}}
As we will examine more closely in Section \ref{sec:limits}, this results from the fact that standard linearized aircraft models include at most one pole at the origin \cite{Cook2007}, \revjh{while the widely used automobile models, which motivated important part of the string stability literature, include two \cite{ploeg2014lp}.}
\revjr{This important difference with automobiles can be attributed to the fact that aircraft actuator dynamics are fast compared to the (relative) motion of the aircraft itself, and can thus be safely neglected in models for timescales relevant to maneuvering in formation.}
However, we find that with a linear state feedback controller such as LQR, velocity disturbances due to wind and wake effects may lead to steady-state errors that degrade the energy savings by pushing the following aircraft away from the optimal positions in the upwash regions of their respective leaders. 
On the other hand, we demonstrate that adding integral control, \rev{a classical approach for eliminating steady-state errors}, can easily result in a formation that is string unstable, \rev{suggesting that string stability may be a challenging design objective when precise relative positioning is important.} 
Finally, we present a control design method
that allows one to achieve both string stability and energy efficiency in the aircraft formation. 

\section{String stability} \label{sec:ss}

\revjh{We now formally define the notion of string stability and discuss some fundamental limitations based on systems theory.} 

\subsection{Definition}

\revjr{Consider a cascaded system of $n$ mobile agents whose positions are denoted by $\p_0(t)$, $\p_1(t)$, $\dots$, $\p_{n-1}(t)$, where $\p_i(t) \in \mathbb{R}^d$ and $d$ is the dimension of the space where the agents move.
Let $\e_i(t) := \p_{i-1}(t) - \p_i(t) - \sepRef$ denote the difference between the actual distance to the preceding agent and a fixed desired separation distance $\sepRef \in \mathbb{R}^d$.
Each agent has a control input $\u_i(t) \in \mathbb{R}^c$, where $c$ is the number of control actuators.
Assume that each agent is modeled by a linear time-invariant (LTI) system such that the dynamics of each agent can be expressed in the Laplace domain as $\p_i(s) = P(s) \u_i(s)$, and that the controller for each agent that follows another depends only on the distance from the immediately preceding agent: $\u_i(s) = C(s) \e_i(s)$.
The separation distance can now be written in terms of the preceding separation distance: $\e_i(s) = T(s) \e_{i-1}(s)$, where $T(s) = (I + P(s)C(s))^{-1}P(s)C(s)$.}
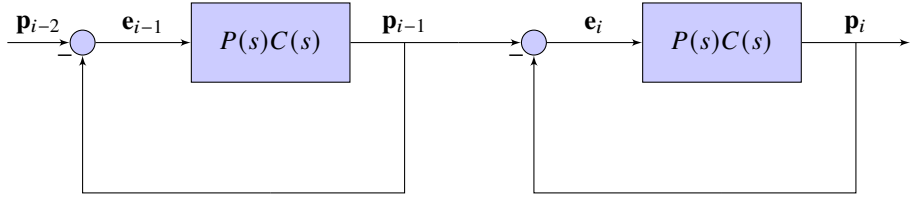
\begin{figure}[ht]
\centering
\begin{tikzpicture}[auto, node distance=2cm,>=latex']
    % Nodes for system 1
    \node [input, name=input] {};
    \node [sum, right of=input] (sum1) {};
    \node [block, right of=sum1, node distance=2.5cm] (sys1) {$P(s)C(s)$};
    \node [output, right of=sys1, node distance=2.5cm] (output1) {};
    \node [feedback, below of=sys1] (feedback1) {};
    
    % Nodes for system 2
    \node [sum, right of=output1] (sum2) {};
    \node [block, right of=sum2, node distance=2.5cm] (sys2) {$P(s)C(s)$};
    \node [output, right of=sys2, node distance=2.5cm] (output2) {};
    \node [feedback, below of=sys2] (feedback2) {};

    % Connections 
    \draw [draw,->] (input) -- node {$\mathbf{p}_{i-2}$} (sum1);
    \draw [->] (sum1) -- node {$\mathbf{e}_{i-1}$} (sys1);
    \draw [-] (sys1) -- node [name=y] {$\mathbf{p}_{i-1}$}(output1);
    \draw [-] (y) |- (feedback1);
    \draw [->] (feedback1) -| node[pos=0.99] {$-$} 
        node [near end] {} (sum1);
    
    \draw [draw,->] (output1) -- node {} (sum2);    
    \draw [draw,->] (sum2) -- node {$\mathbf{e}_i$} (sys2);
    \draw [->] (sys2) -- node [name=y] {$\mathbf{p}_i$}(output2);
    \draw [-] (y) |- (feedback2);
    \draw [->] (feedback2) -| node[pos=0.99] {$-$} 
        node [near end] {} (sum2);
\end{tikzpicture}
\caption{Leading and following agents as cascaded LTI system.} \label{fig:sensitivity}
\end{figure}
\revjr{To see why this is true, observe in Fig. \ref{fig:sensitivity} that $\e_i(s) = \p_{i-1}(s) - \p_i(s)$ and \rev{$\p_i(s) = P(s)C(s)\e_i(s)$}, resulting in
\begin{align*}
    \e_i(s) = P(s)C(s)(\e_{i-1}(s) - \e_i(s)) \implies \e_i(s) =  (I + P(s)C(s))^{-1}P(s)C(s) \e_{i-1}(s) = T(s) \e_{i-1}(s).
\end{align*}}

\rev{The cascaded system is said to be \emph{string stable} if $P(s)C(s)$ is stable and if \revjh{no frequency of disturbance is amplified from leader to follower, i.e.} 
$\sup_\omega |T(j\omega)| \leq 1$, for single-input single-output \revjr{(SISO)} systems \cite{ploeg2014lp}.} 
For multiple-input multiple-output systems (i.e., $T(s)$ is a transfer matrix), the generalized string stability criterion is $\sup_\omega \bar\sigma[T(j\omega)] \leq 1$, where $\bar\sigma$ denotes the maximum singular value.
\revTwo{In systems that are not string stable, there is thus a frequency that gets amplified at each following agent, leading to an exponential growth along the cascaded system of disturbances at that frequency.}

\subsection{Known limitations} \label{sec:limits}

In classical feedback control theory, a phenomenon known as the \textit{waterbed effect} places fundamental theoretical limitations on the ability to simultaneously achieve good tracking performance while attenuating disturbances across the frequency spectrum.
This comes as a direct consequence of Bode's integral constraint on the sensitivity function.
Let $G(s):=P(s)C(s)$ be the open loop transfer function of a controller $C(s)$ applied to the SISO system $P(s)$, where $s\in\mathbb{C}$.
The sensitivity function measures the effect of the reference signal and input disturbances on the error signal and is given by $S(s) = \frac{1}{1 + G(s)}$.
If $G(s)$ has at least two more poles than zeros and no poles in the right half-plane, then Bode's constraint can be expressed as follows:
\begin{align*}
    \int_0^\infty \ln{|S(j\omega)|}d\omega = 0.
\end{align*}

A related yet somewhat lesser known result involves the complementary sensitivity function $T(s) = \frac{G(s)}{1 + G(s)}$, which turns out to be highly relevant for string stability, since the amplification of disturbances to the error signal from leader to follower is given exactly by this function $T(s)$.

The following definition uses the notion of \textit{system type}, which is the number of pure integrators in the open-loop transfer function $G(s)$.
That is, a system written in the form 
\[G(s) = \frac{1}{s^\ell}\frac{K(T_1 s+1)(T_2 s+1)...(T_{n_z} s + 1)}{(T_a s+1)(T_b s+1)...(T_{n_p} s + 1)}
\]
is said to be of type $\ell$.
We will also need to define the \textit{velocity error constant}, which is \rev{the inverse of} the steady-state error of a system in response to a unit ramp input, and is given by $K_v := \lim_{s\rightarrow 0} sG(s)$.
The integral constraint on the complementary sensitivity function is then given by
\begin{align} \label{eq:sT}
    \int_0^\infty \frac{\ln{|T(j\omega)|}}{\omega^2}d\omega 
    = -\frac{\pi}{2 K_v} + \pi\sum_{i=1}^{q} \frac{1}{z_i},
\end{align}
where $z_1,\dots,z_q$ denote any zeros of $G(s)$ in the open right half of the complex plane.
For systems of type 0, the right side of \eqref{eq:sT} is infinity.
For systems of type 2 \revTwo{and higher} with no zeros in the right-half of the complex plane,
it is equal to zero \rev{\cite{emami2019bode}}.

It was shown in \cite{seiler2004disturbance} that no linear controller can render a cascaded formation string stable if all agents are identical, LTI, SISO, strictly proper \revTwo{(i.e. the degree of the denominator is higher than the degree of the numerator)}, have two poles at the origin, and only measure the distance to the preceding agent (i.e. relative velocity is not available).
\rev{This result holds because in order for the integral in \eqref{eq:sT} to be equal to zero, $\frac{|T(j\omega)|}{\omega^2}$ must be uniformly equal to one, or else must be greater than one for some frequencies, which translates to string instability at those frequencies.}
This result was extended in \cite{middleton2010string} to show that system heterogeneity and an extended but limited amount of forward communication are not sufficient to avoid string instability, although they may improve performance.
In the next section, we briefly describe some methods for modifying the system in a way that makes string stability achievable when the fundamental limitation is in effect.

\subsection{Known remedies} \label{sec:remedies}

Several methods have been proposed to \rev{mitigate string instability, for example by also using the distance to the following vehicle in the controller \cite{lestas2007scalability,middleton2010string} and by using heterogeneous controllers \cite{barooah2007control}. 
Other approaches can avoid the problem altogether, including the use of a} sufficiently large time headway (controlling inter-agent time in addition to inter-agent distance) \cite{yanakiev1996simplified}, knowledge of the control input of the preceding agent \cite{ploeg2014lp}, and knowledge of the position of the absolute leader \cite{swaroop1997string}.

\rev{From the known methods for avoiding string stability, we focus here on the use of a time headway, since that requires no inter-agent communication,} \revTwo{which would significantly increase the complexity.}
The term \textit{time headway} here refers to a reference separation distance that is not constant but rather depends on the agent's own velocity.
This can be expressed in general terms as $\sepRefElem(t) := \sepRefElem_0 + \headway \velElem_i(t)$, where $\sepRefElem(t)$ denotes the reference separation, $\sepRefElem_0$ denotes a fixed distance, $\headway$ denotes the time headway constant, and $\velElem_i(t)$ denotes the velocity of agent $i$ at time $t$.
It is well-established that using such a reference separation with a sufficiently large value of $\headway$ can lead to string stability of a formation \cite{klinge2009string}.
We will return to this result in Section \ref{sec:aircraft}, but first we introduce the model we will use to examine string stability in the context of aircraft formations.

\section{Aircraft model} \label{model}

The simulated aircraft is based on the Airbus A320 flying at 230 m/s, with a wingspan of $b=$ 34 m and other dimensions listed in Table \ref{tab:aircraftDimensions}.
\rev{This aircraft was chosen because its typical medium-distance flights provide it the opportunity to join formations which may last for hours, potentially yielding substantial energy savings. The A320 also uses the typical fixed-wing design, thus the results found with this aircraft should generalize well.}
The aircraft's controls are the deflections of its rudder (vertical tail), elevator (horizontal tail), and oppositely-moving ailerons (outer wing flaps), as well as the engine thrust. This section describes the simulation of the A320 and its wake, the perturbations it encounters, and the drag reduction which is possible in formation flight. 

\subsection{Vortex Lattice Method}
The forces and moments on the aircraft are computed with a steady Vortex Lattice Method (VLM) \cite{katz_plotkin}, which omits compressibility effects.  This approach uses potential flow theory to model the aircraft as a collection of vortex line filaments. Each of these induces the flow to spin around the filament, as defined in Equation \eqref{eq:vortexFilament}.
The strength of these filaments is determined by solving a linear set of equations that ensures that flow is tangent to the surfaces at certain control points. The forces on each vortex filament on the aircraft, which are a product of the external flow velocity and the filament's circulation, are compiled to yield the total aerodynamic forces and moments on the aircraft. The external velocity field includes the flow due to motion of the aircraft, the wakes of other aircraft, constant wind, or unsteady turbulence. Although real aircraft would not fly in formation with a constant crosswind, it is the limiting case of a very long-wavelength turbulent gust, and its steadiness enables us to see important effects of such a gust. \rev{As this is an inviscid potential flow model, it does not estimate forces like the parasitic drag. This is captured by the drag coefficient defined in Table \ref{tab:aircraftDimensions}.}

\begin{figure}[ht]
	\centering
	\includegraphics[width=0.7\linewidth]{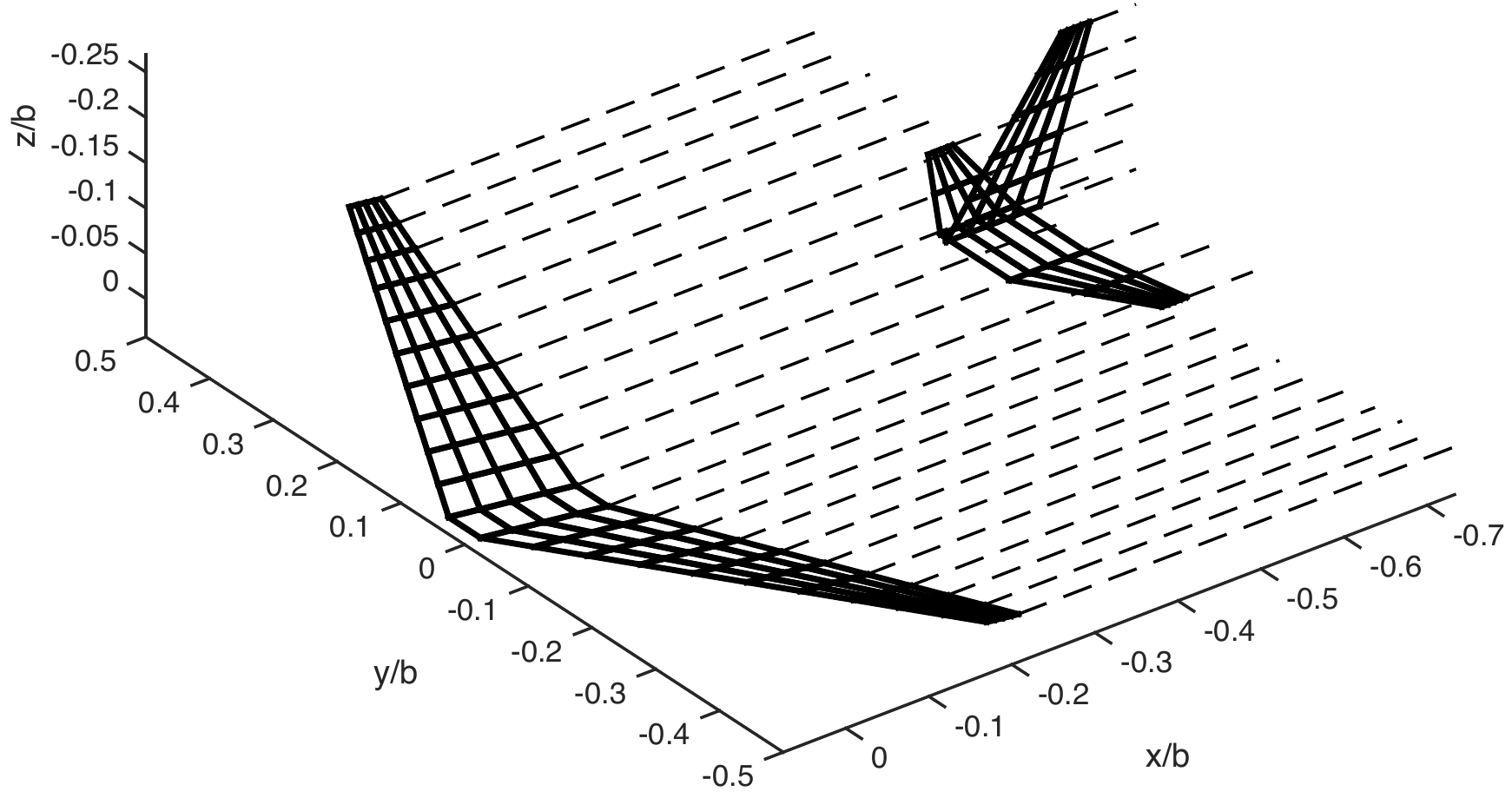}
	\caption{VLM model of an Airbus A320. The dashed lines represent the trailing vortices, and the solid lines are the edges of the vortex panels.}
	\label{fig:airplaneVLM}
\end{figure}

The wing and tail surfaces are represented by vortex panels composed of four vortex segments. Each panel at the trailing edge of the surfaces also has a horseshoe vortex. These are made of one finite segment at the trailing edge, and two semi-infinite segments which extend in the $-x$ direction. The effects of control surface deflections are approximated by tilting the surface normal vectors of the appropriate vortex panels. The vortex filaments are shown in Fig. \ref{fig:airplaneVLM}, along with the axis conventions. The streamwise direction $x$ is positive forward, the lateral direction $y$ is positive to the aircraft's right side, and the vertical direction $z$ is positive down.

\subsection{Wake velocity}
An aircraft in flight leaves behind a region of disturbed air known as the wake. Although the wake is complex near the aircraft, after a few wingspans downstream the wake `rolls up' into a pair of counter-rotating vortex lines, as seen in Fig. \ref{fig:wakeVelocity}. The flow due to the rolled-up vortices is well approximated by a horseshoe vortex model \cite{widnall1975} with a velocity field as in Equation \ref{eq:wakeVel}. The distance between the wake vortices is a fraction of the wingspan: $b \frac{\pi}{4}\approx 0.79b$. The velocity field of the simplest horseshoe vortex model is singular on the vortex lines, but Equation \ref{eq:wakeVel} includes a core size which removes the singularity. The maximum speed is found near that core radius, which is chosen as 5\% of the wingspan.
This horseshoe vortex is used to model the influence of a leading aircraft on its followers.
Because the aircraft maneuvers are relatively slow in comparison to the cruise velocity, and the weight of the aircraft is approximated as constant, the strength of the wake vortices are constant over time. 

In reality, the wake would not be perfectly straight: it has a varying position that depends on the motions of the aircraft. With a streamwise separation of 10 wingspans, such disturbances would take 1.48 seconds to propagate to the next aircraft in line. We model this propagation effect as a delay in the position of the horseshoe vortex. Thus, the relative position of the leader 1.48 seconds ago determines the current velocity field at the follower. Because the relative changes in $x$-position are small, we assume that this delay has a constant value.

\subsection{Wake exploitation}
\label{sec:wakeSpacing}
If an aircraft flies in a region of upwash, it can reduce its drag substantially. As a demonstration of the effect, we use steady lifting line theory \cite{katz_plotkin} to approximate the drag coefficient on an elliptic wing in steady level flight as
\begin{align}
    C_D = C_{D,0} + \frac{C_L^2}{\pi AR} - C_L \frac{w}{U},
\end{align}
where $C_{D,0}$ is the zero-lift drag coefficient, $C_L$ is the lift coefficient, $AR$ is the aspect ratio of the wing, $w$ is the uniform upwash on the wing, and $U$ is the speed of the aircraft. This shows that maximizing the upwash on the wing minimizes its drag.

When flying behind another aircraft, the optimal position for an aircraft to maximize the upwash is to fly behind and to the left or right of the leader, with its wingtip touching one of the wake vortices. That is, at a lateral separation of $b(1+ \pi/4)/2 \approx 0.89 b$, and at the same altitude as the wake. 

\subsection{Turbulence}
The simulations allow for in-flight turbulence, modeled as von Karman Turbulence \cite{matlabVK}. This is a stochastic method that uses white noise to generate spectrally accurate turbulence along the $x$ axis. The turbulence is assumed to be `frozen,' and so is fixed in space. The aircraft's surfaces experience the velocity fluctuations as a function of their $x$-positions, so the wings feel a gust before the tail does. The generated turbulence has a length scale of 762 m (22.3$b$) and an intensity which is set for each simulation. Because the relative $y$ and $z$ separations of the aircraft are on the scale of one wingspan, it is appropriate to use this one-dimensional model of turbulence.

\subsection{Aircraft dynamics}
For the aircraft flying in the absence of external disturbances, the system dynamics (linearized around the trimmed state \cite{leng}) are
\begin{equation}
\dot{\mathbf{x}} \approx A \mathbf{x} + B \mathbf{u},
\end{equation}
where $\mathbf{x}$ is the state of the aircraft, and $\mathbf{u}$ is the control input. \revThree{The state $\mathbf{x} := [\p^\top\, \v^\top\, \a^\top\, \dot\a^\top]^\top$ is composed of the aircraft's three-dimensional position $\p = [x\,\, y\,\, z]^\top$ relative to some fixed global reference point, three-dimensional Euler angles $\a = [\phi\,\, \theta\,\, \psi]^\top$ that define its orientation, and its rates of translation $\v$ and rotation $\dot\a$, with the trimmed values subtracted from each component.} The control vector $\mathbf{u}$ has the change in thrust, and the deflections of the ailerons, rudder, and elevators. The linearized dynamics matrices, $A$ and $B$, were calculated from the VLM model using central differences. The linearized dynamics matrices, the state vector, and the control vector are presented in Appendix \ref{app:systemDynamics}. 

The behavior of simple fixed-wing aircraft can be separated into longitudinal and lateral dynamics, which are linearly uncoupled \cite{Cook2007}, simplifying the control problem. The longitudinal dynamics involve the $x$ and $z$ positions and velocities, as well as the pitch and pitch rate. The lateral dynamics involve the $y$ position and velocity, as well as the yaw and roll angles and rates. 

Similarly to Binetti et al.\cite{binetti2003formation}, the effect of the external velocity field is treated as a nonlinear exogenous input, $\mathbf{w}$. For ease of computation, the effects of the external velocity are computed using the VLM with the aircraft in the disturbance-free trimmed state. %This yields $\mathbf{w}=\dot{\mathbf{x}}_{ext.\, vel}-\dot{\mathbf{x}}_{trim}$, which 
We include this in the state evolution equation as:
\begin{equation} \label{eq:dyn}
\dot{\mathbf{x}} = A \mathbf{x} + B \mathbf{u} + \mathbf{w}.
\end{equation}

For simplicity, we assume that the aircraft have full knowledge of their state. This assumption is reasonable with modern instrumentation; measurements of the accelerations and angles are available with high precision with onboard gyroscopes and accelerometers, and the relative positions of the aircraft can also be measured precisely.

\section{String stability and energy savings in aircraft formations} \label{sec:aircraft}

\rev{Based on the model developed in the previous section, we proceed to analyze the string stability and energy savings of a group of airplanes in a diagonal line formation.} 

\subsection{Problem formulation}
Consider a formation of $n$ airplanes 
and let the dynamics of each airplane $i \in\{1,\dots,n\}$ be modeled by \eqref{eq:dyn}.
Suppose that the desired location with respect to a leading aircraft is chosen to be a constant reference vector $\sepRef := [\sepRefElem_x,\sepRefElem_y,\sepRefElem_z]^\top$ and assume that the target formation is a sequence of aircraft who maintain this reference separation with respect to the preceding aircraft.
Specifically, each aircraft's target position is offset from its immediate leader by $\sepRefElem_x=10b$ (10 wingspans) in the downstream direction, $\sepRefElem_y=0.89b$ laterally, and with the same altitude ($\sepRefElem_z=0$). 
The streamwise separation was chosen to be large enough to avoid the risk of collisions, while reducing the magnitude of wake meandering. This optimal lateral separation distance is derived in section \ref{sec:wakeSpacing}. 
Each plane experiences the wake of its immediate leader with a delay of 1.48s\rev{, due to the formation's speed and streamwise separations}. Ideally, each plane would hold its position perfectly with respect to its leader and so enjoy substantial drag reduction. Toward that goal, the aircraft need appropriate control methods.

Given the target offset vector $\sepRef$, the error of each airplane with respect to the this reference separation can be written as $\err_i := \pos_{i-1} - \pos_i - \sepRef$.
Note that since trimmed (steady-state) values are subtracted from the state, the aircraft velocity state represents the deviation $\v = \dot \p - \v_0$ from the trimmed velocity vector $\v_0$. The control objective is then to stabilize each separation error $\err_i$ to near zero while achieving string stability of the formation, as defined in Section \ref{sec:ss}.

\subsection{Comparison to the vehicle platooning problem}

This problem resembles that of vehicle platooning\rev{, i.e. controlling a sequence of cars to follow each other at close distance to increase traffic throughput while saving fuel,} but in three dimensions instead of one.
One of the reasons that string stability is \rev{widely considered to be a challenge in the context of} vehicle platoons is that \rev{after applying nonlinear  control to linearize the automobile dynamics in a standard technique known as \textit{feedback linearization}, the resulting models} typically include two pure integrators \cite{swaroop1994comparison}.
This means that there is no linear controller that can achieve string stability of the vehicle platoon under the conditions described in Section \ref{sec:limits}.
In contrast to vehicle platoons, the linearized aircraft dynamics considered here include only one pure integrator (in each spatial dimension), due to the fact that the aircraft actuator dynamics are assumed to be fast \revjr{enough to be neglected in the context of the dynamics of the aircraft as a whole.}
\rev{This means that aircraft models based on these dynamics may indeed be rendered string stable by a linear controller.}

For example, suppose we wish to achieve string stability of the lateral ($y$) dynamics in the aircraft formation.
\revThree{Aircraft dynamic models are often separated into lateral and longitudinal components, since these two subsystems are generally independent after linearization. We can therefore consider the lateral subsystem in isolation for the purposes of this analysis.}
Using the model from Section \ref{model}, the transfer functions from the ailerons and rudder to the $y$ position are of system type one, meaning they each include one pure integrator (see Appendix \ref{app:transferFunctions}).
Therefore, the limitations described in Section \ref{sec:limits} do not apply, and it may be possible to design a string stabilizing controller for the relative lateral positions of such an aircraft formation.

\subsection{A string stabilizing state feedback controller}

\begin{figure}[ht]
	\centering
	\begin{subfigure}[t]{.4\linewidth}
    \centering
     \includegraphics[width=\linewidth]{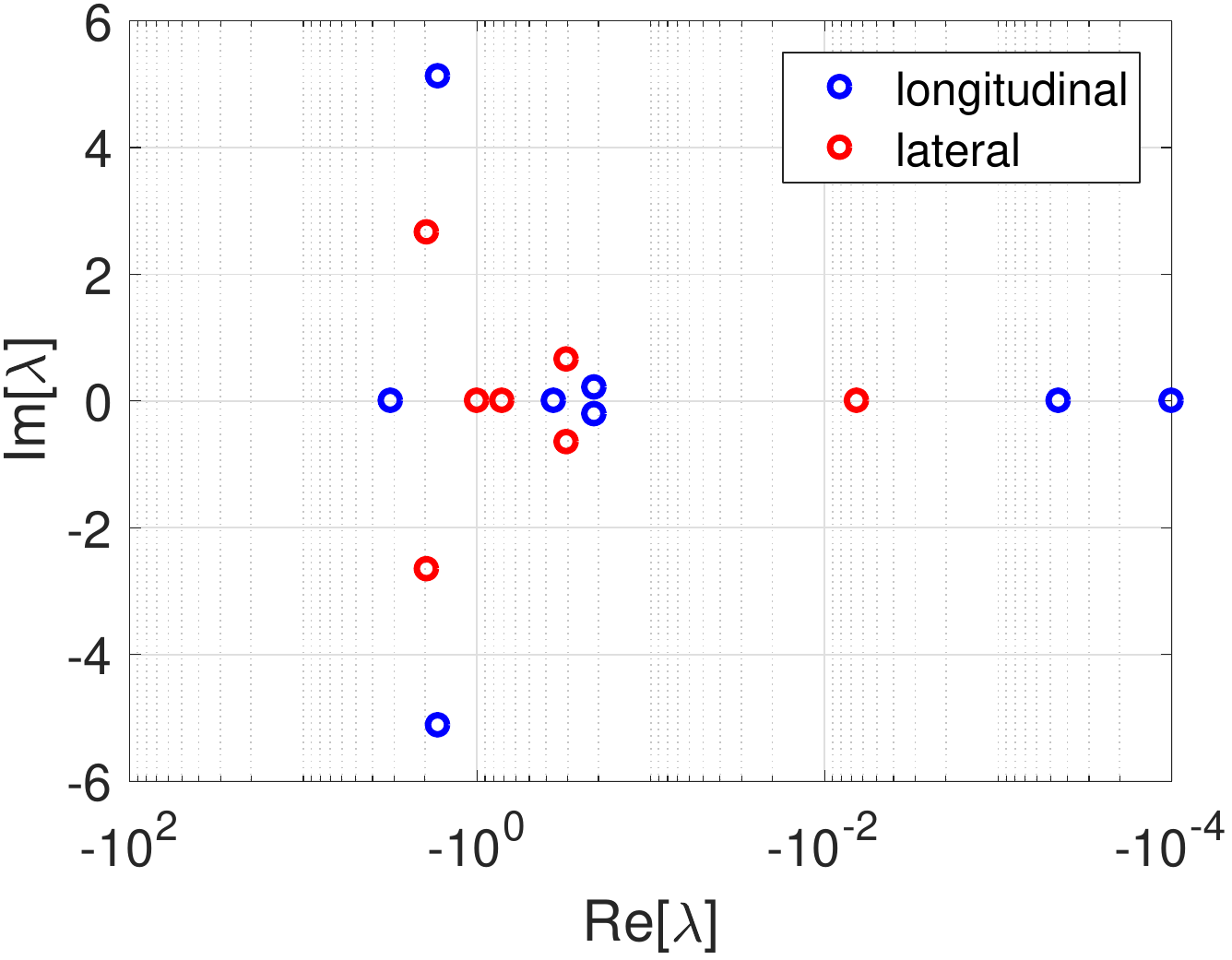}
		\caption{\small{Eigenvalues of closed-loop system}}\label{fig:eigbodeLQR-a}		
	\end{subfigure}
	\quad
	\begin{subfigure}[t]{.55\linewidth}
    \centering
     \includegraphics[width=\linewidth]{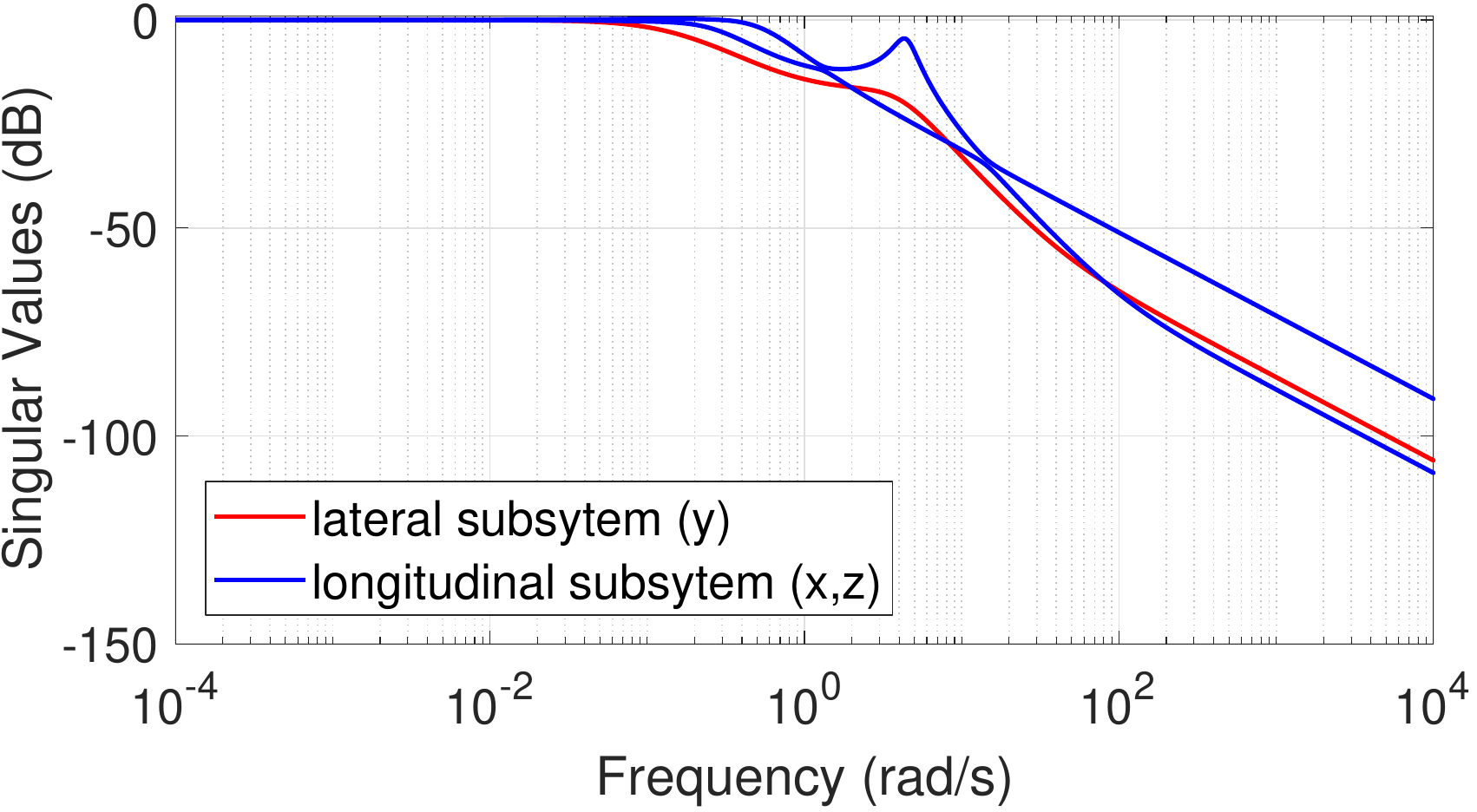}
		\caption{\small{Singular values of $T(s)$}}\label{fig:eigbodeLQR-b}
	\end{subfigure}
    \caption{Closed-loop eigenvalues and singular values of complementary sensitivity transfer matrix for LQR controller.}
    \label{fig:eigbodeLQR}
\end{figure}

\rev{We show that it is indeed possible to achieve string stability of the formation in all three dimensions by designing a standard linear quadratic regulator (LQR) for optimal disturbance attenuation from leading to following aircraft.
The stability of the closed-loop system and the string stability of the formation are illustrated in Figure \ref{fig:eigbodeLQR}. The singular values of the complementary sensitivity transfer matrix $T(j\omega)$ do not exceed one, and therefore the formation is string stable (See Appendix \ref{app:controlSpecs} for the transfer functions and control gains used in this example).} 

\rev{We have thus designed an LQR controller that renders the aircraft formation string stable, but} another important factor to consider is the nature and impact of disturbances on the energy-saving performance of aircraft formations.
\rev{While disturbances due to wind are often neglected in vehicle platoon models, it is not justifiable to neglect such disturbances for aircraft, particularly for the level of precision required to achieve good energy savings, as we will see in the next section.}

\subsection{Trade-off between string stability and energy savings} \label{sec: tradeoff}
Since the open-loop system with state feedback in the previous section includes one integrator, a constant disturbance will result in zero steady-state error.
However, wake effects and wind gusts are disturbances that impact the velocity of an aircraft. For example, and a constant wind disturbance could equivalently be thought of as a ramp disturbance to the position of the aircraft.
Linear systems theory dictates that in order for a closed-loop system to have zero steady-state error in the presence of a ramp input, its system type needs to be at least two, i.e., the open loop transfer function should contain at least two pure integrators \cite{seiler2004disturbance}.
Otherwise, wind disturbances will result in a degradation of tracking performance, and the same holds true for wake effects from preceding aircraft. 

\begin{figure}[ht]
	\centering
	\begin{subfigure}[t]{.45\linewidth}
    \centering
    \includegraphics[width=\linewidth]{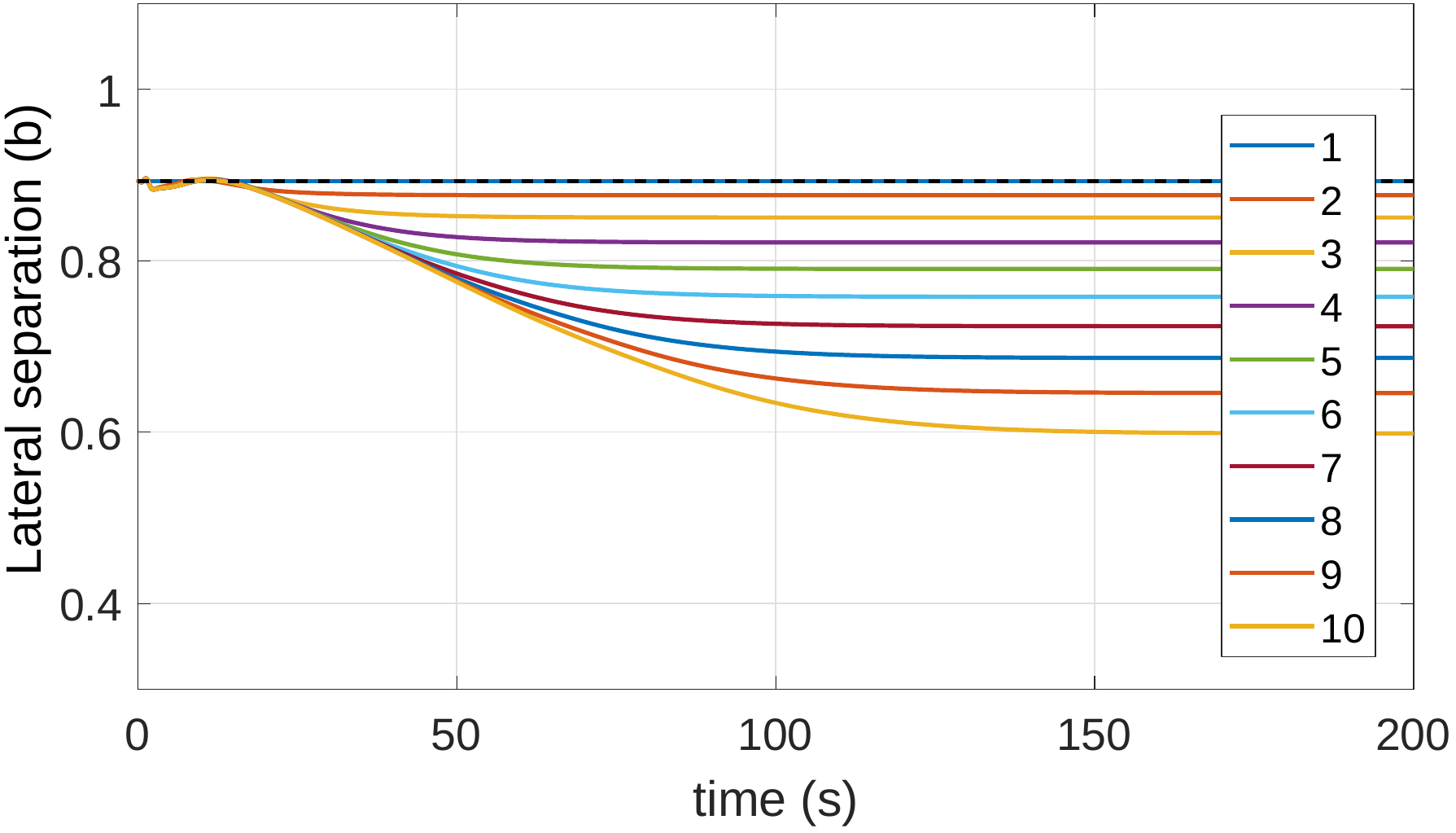}
		\caption{\small{Relative lateral positions}}\label{fig:wakeLqr-a}		
	\end{subfigure}
	\quad
	\begin{subfigure}[t]{.45\linewidth}
    \centering
    \includegraphics[width=\linewidth]{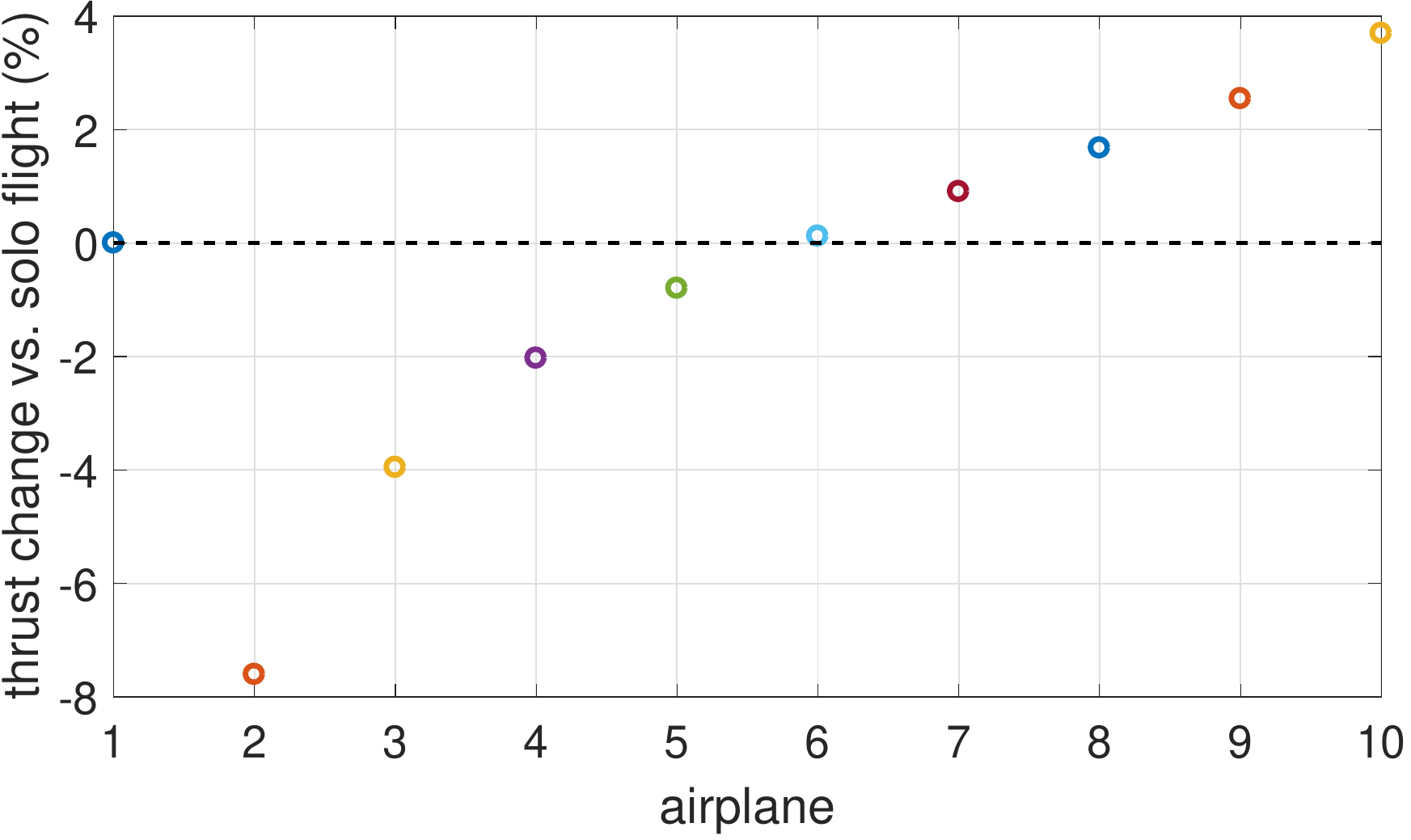}
		\caption{\small{Final change in thrust}}\label{fig:wakeLqr-b}
	\end{subfigure}
    \caption{Simulations of the LQR controller show that disturbances due to wakes can cause increasing steady-state errors along the aircraft formation. This degrades the energy savings, which is measured as the percentage change in thrust with respect to solo flight.%
    }
	\label{fig:wakeLqr}
\end{figure}

Fig. \ref{fig:wakeLqr-a} shows that significant steady-state errors emerge in the lateral positions of the 10-aircraft formation using the LQR control designed above when wake effects are included in the model. This simulation was initialized with the aircraft in their ideal formation, but the presence of the wakes pushed them out of position. The resulting final positions are thus dependent on the initial conditions and subsequent trajectories.
While a position tracking error of this magnitude might be acceptable for solo flight, in formation flight the resulting misalignment with the upwash region of the preceding aircraft's wake significantly degrades the energy savings. In this case, each successive aircraft drifts further into the downwash region of its respective leading aircraft, resulting in the increasing thrust shown in Fig. \ref{fig:wakeLqr-b}. 

\begin{figure}[ht]
	\centering
	\begin{subfigure}[t]{.45\linewidth}
    \centering
    \includegraphics[width=\linewidth]{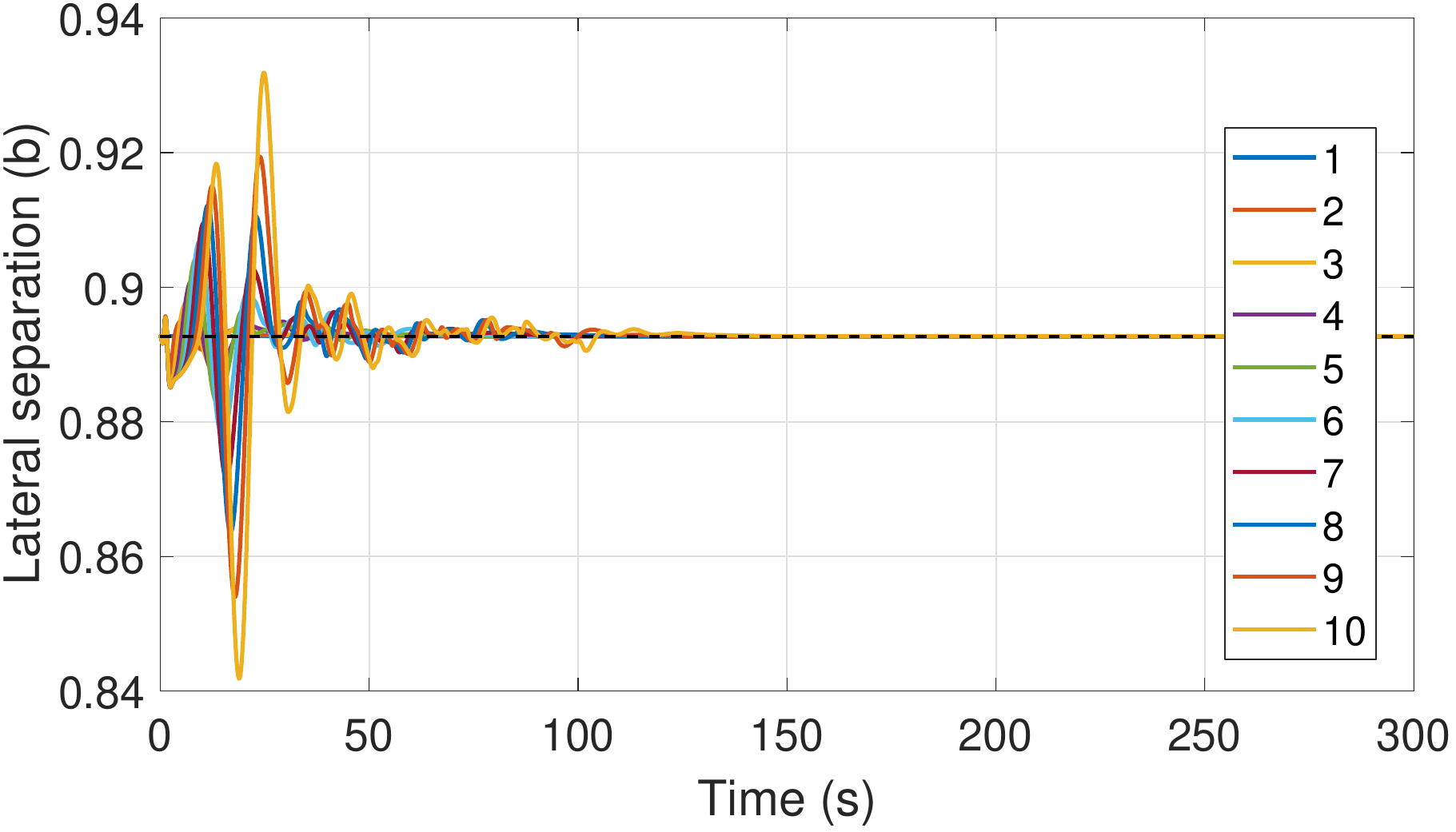}
		\caption{\small{Relative lateral positions}}\label{fig:wakeLqrInt-a}		
	\end{subfigure}
	\quad
	\begin{subfigure}[t]{.45\linewidth}
    \centering
    \includegraphics[width=\linewidth]{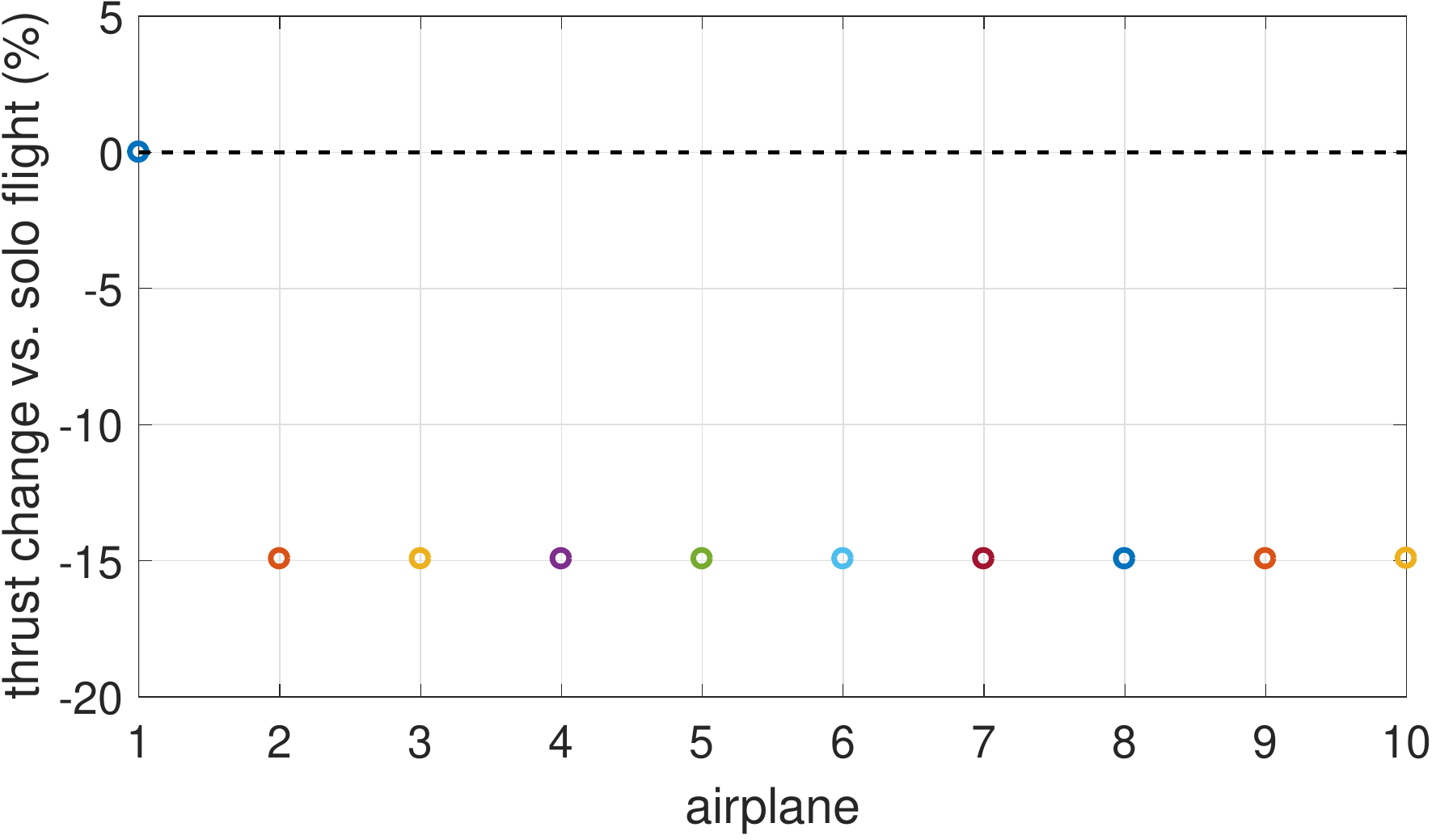}
		\caption{\small{Final change in thrust vs. solo flight}}\label{fig:wakeLqrInt-b}
	\end{subfigure}
    \caption{The LQR plus integral controller eliminates the steady-state error, yielding constant energy savings for all following aircraft. However, some amplification in the transient response is observed, which suggests that the formation is string unstable.}
	\label{fig:wakeLqrInt}
\end{figure}

The standard way to resolve this is to add integral control action, which can be achieved by integrating the relative position in each spatial dimension and including these additional three states in the LQR design. 
The performance of the LQR plus \rev{integral} controller in the presence of wake effects is shown in Fig. \ref{fig:wakeLqrInt}. After the initial transients, the followers enjoy uniform and significant drag reduction.

One can already observe a potential problem with this design, however, which is that there is some overshoot in response to the wake disturbances that is amplified from one airplane to its follower, suggesting the presence of string instability. Indeed, the result of adding integral control is an open-loop transfer function that has two integrators in each spatial dimension (see Appendix \ref{app:transferFunctions}).
We see in Fig. \ref{fig:bodemag_lqrint} that the magnitude of the diagonals of the complementary sensitivity transfer matrix for the chosen LQR plus integral controller exceed one in each dimension.
\begin{figure}[ht]
    \centering
    \includegraphics[width=0.5\linewidth]{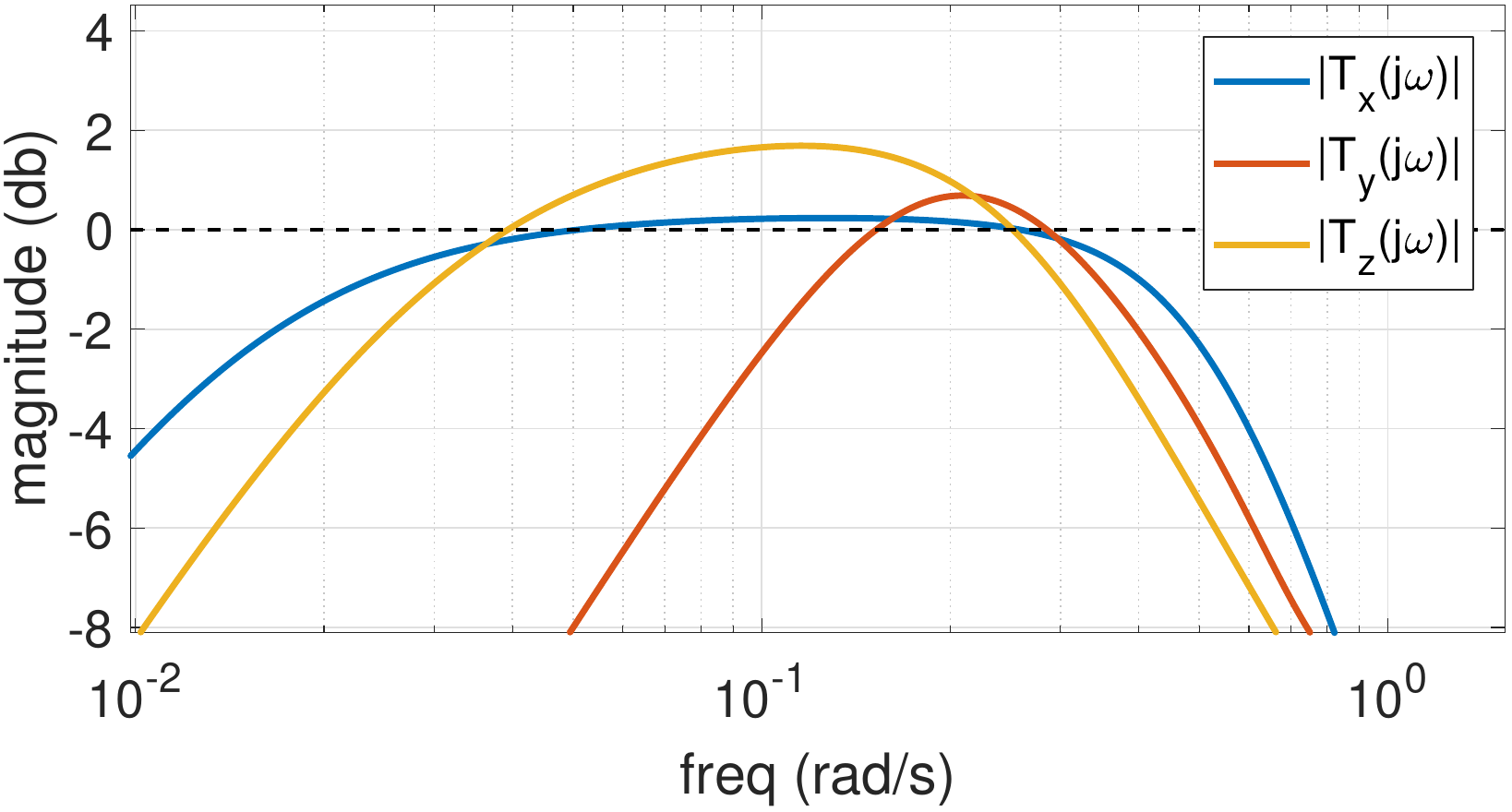}
    \caption{For the LQR plus integral controller, the magnitude of the diagonals of the complementary sensitivity transfer matrix $T(s)$ exceed one in each dimension near 0.1 rad/sec, which indicates that the formation is string unstable.}
    \label{fig:bodemag_lqrint}
\end{figure}
\begin{figure}[ht]
    \centering
    \includegraphics[width=0.6\linewidth]{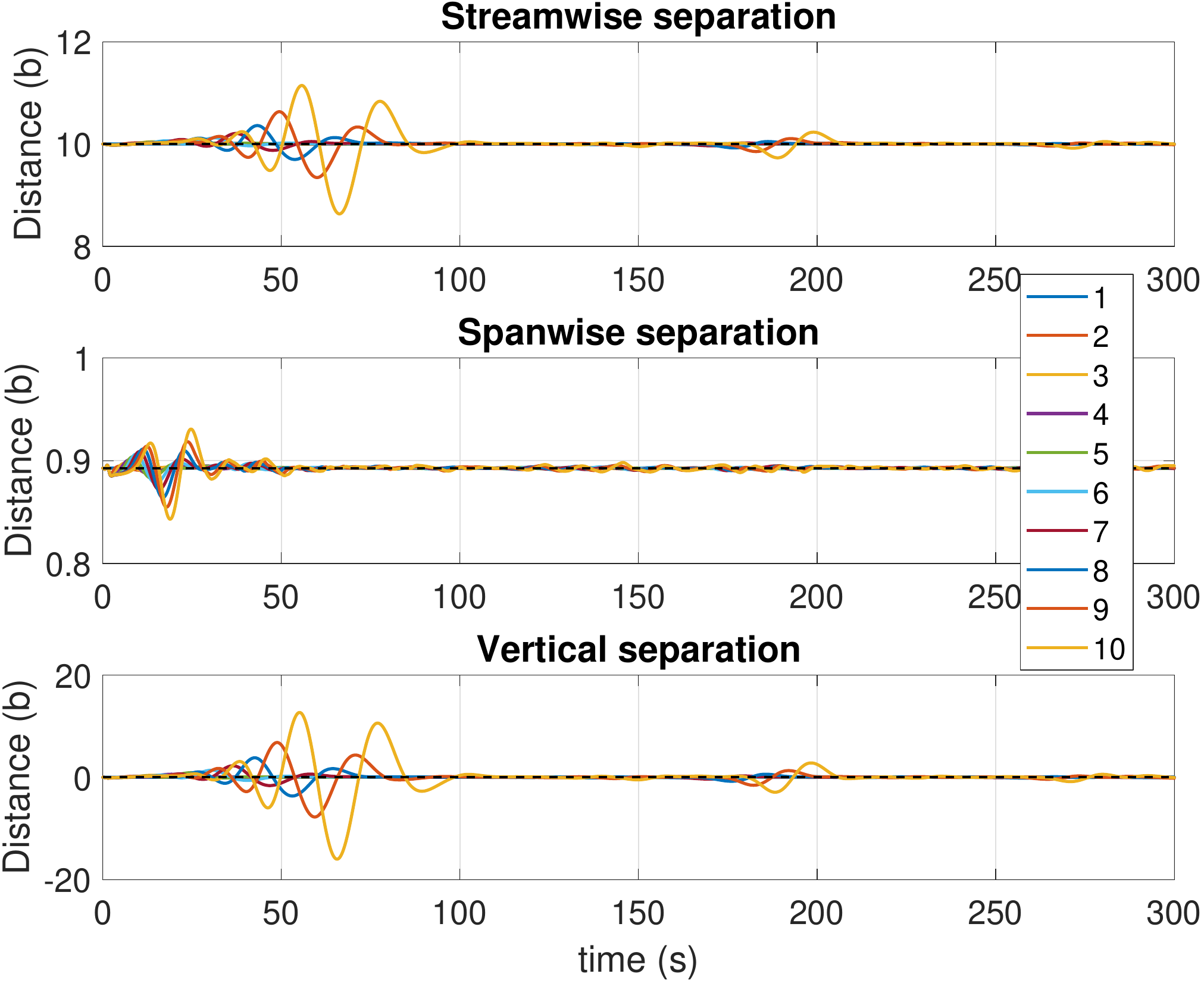}
    \caption{Using the LQR plus integral controller, a formation of 10 planes flying in 2\% turbulence intensity exhibits \revjh{occasional} large deviations from the prescribed relative positions.}
    \label{fig:turb_lqrInt}
\end{figure}
\begin{figure}[ht]
    \centering
    \includegraphics[width=0.5\linewidth]{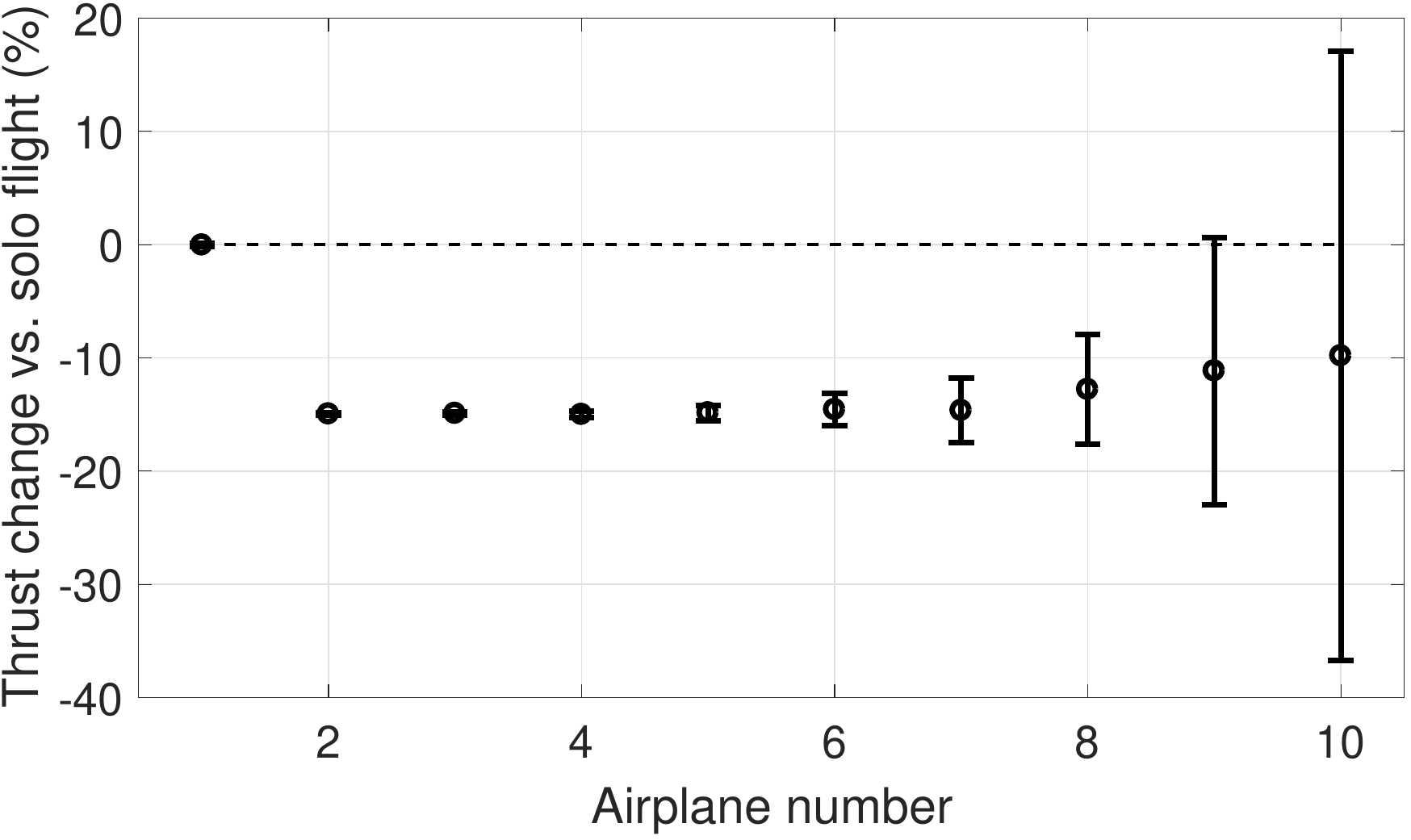}
    \caption{Deviations from desired offset positions due to string instability \revThree{in 2\% turbulence} degrade the energy-saving performance of the formation using the LQR plus integral controller. The error bars show the standard deviation in the thrust change vs. solo flight.}
    \label{fig:turb_lqrInt_energy}
\end{figure}

Since the system is stable in the classical sense, the formation eventually recovers from the initial transients and converges to the prescribed formation.
However, a persistent disturbance such as ambient turbulence may lead to more serious problems.
Fig. \ref{fig:turb_lqrInt} shows the relative positions of 10 airplanes in formation with LQR plus integral control in the presence of 2\% turbulence intensity. The stochastic nature of turbulence occasionally excites the string instability, resulting in large deviations from the optimal positions that \rev{degrade the energy-savings of the formation, as shown in Fig. \ref{fig:turb_lqrInt_energy}. 
In addition, such oscillations correspond to undesirable and perhaps unsafe flying conditions.}

In summary, we find that string stability may be challenging to attain in aircraft formations whose accurate relative positioning is important, for example when the objective is to fly in the upwash of preceding agents. 
\revjr{Since our addition of integral control resulted in an open-loop transfer function of controller plus aircraft dynamics that has two pure integrators in each spatial dimension, the fundamental limitations described in Section \ref{sec:limits} are now in effect.
Indeed, since the dynamics need two integrators for accurate relative position tracking in the presence of wind disturbances, the string stability limitation holds for any controller that achieves this objective.}
However, recall that one of the remedies for string instability discussed in Section \ref{sec:remedies} was to add a time headway $\delta(t) := \delta_{ref} + \headway \velElem_i(t)$ to the reference position.
Standard state feedback already includes a velocity feedback term in each dimension, which can be used to introduce a time headway in the final controller.
\revjr{It was shown in \cite{klinge2009time}
that a sufficiently large time headway can render a cascaded formation string stable even when the component subsystems have two integrators.
Moreover, since the desired steady-state is a formation flying at a constant velocity, the 3-dimensional reference separation with a time headway term $\boldsymbol{\delta}(t) := \boldsymbol{\delta}_{ref} + \headway \vel_i(t)$ will converge to the desired constant separation distance $\boldsymbol{\delta}_0$ as the trimmed velocity $\vel_i(t)$ goes to zero.
This leaves open the possibility that both string stability and accurate tracking performance are achievable in aircraft formations.
In the next section, we show one approach for designing a controller to achieve both of these objectives.}

\rev{\section{Design of a string stabilizing controller with good energy saving performance}}

So far we have seen that integral control is necessary for accurate relative position tracking in the presence of aerodynamic disturbances, but that the resulting system may become string unstable.
In this section, we will show how to design for both string stability and tracking performance by tuning the control gains \revThree{(}including the integral and velocity feedback terms\revThree{)} such that the complementary sensitivity function satisfies the string stability constraint, while also ensuring zero steady-state error.
Since both integral control and velocity feedback are present in the formation control design of \cite{binetti2003formation}, which showed good qualitative performance, we adopt a similar architecture and \revThree{use structured H-infinity synthesis} to achieve our control objectives.

\begin{figure}[ht]
\centering
\tikzstyle{block} = [draw, rectangle, 
    minimum height=2em, minimum width=1em]
\tikzstyle{sum} = [draw, circle, node distance=1cm]
\tikzstyle{input} = [coordinate]
\tikzstyle{output} = [coordinate]
\tikzstyle{pinstyle} = [pin edge={to-,thin,black}]

\begin{tikzpicture}[auto, node distance=1.2cm,>=latex']
	% We start by placing the blocks
	\node [input, name=input] {};
	%
	% separation tracking loop
	\node [sum, right of=input] (sumOuter) {};
	\draw [->] (input) -- node[pos=0.99] {$+$} node[pos=-.05] {$\mathbf{p}_{i-1}$} (sumOuter);
	\node [block, above=5mm of sumOuter] (delta) {$\boldsymbol{\delta}_\text{ref}$};
	\draw [->] (delta) -- node[pos=0.8] {$-$} (sumOuter);
	\node [block, right=10mm of sumOuter]  (KP) {$\frac{1}{s}K_v K_p + K_v K_d$};
	\draw [->] (sumOuter) -- node[pos=0.6] {$\mathbf{e}_i$} (KP);
	%
	% middle loop
	\node [sum, right=5mm of KP] (sumMiddleLoop) {};
	\draw [->] (KP) --  node[pos=0.5] {$\tilde{\mathbf{u}}_i$} node[pos=0.99] {$+$} (sumMiddleLoop);
	%
	% inner loop
	\node [sum, right=10mm of sumMiddleLoop] (sumInner) {};
	\draw [draw,->] (sumMiddleLoop) -- node[pos=0.99] {$+$}  (sumInner);
	\node [draw, rectangle, minimum height=4em, minimum width=1em, right of=sumInner,node distance = 1.5cm] (aircraft) {Aircraft};
	\draw [->] (sumInner) --node {$\mathbf{u}_i$} (aircraft.west);
	\node [block, below of=aircraft] (Kstar) {$K_{\alpha}$};
	\draw [->] (aircraft.-30) -|  ++(0.5,0) |- node[pos=0.7] {$\bar\a_i$}  (Kstar.east);
	\draw [->] (Kstar.west) -| node[pos=0.99] {$-$} (sumInner);
	%
	% aerodynamic interference
	\node [block, above of=aircraft] (AeroInterfere) {Wakes};
	\draw [->] (aircraft.35) -| ++(0.5,0) |- (AeroInterfere.east);
	\draw [->] (AeroInterfere.west) -| ++(-0.5,0)  node[pos=0.8] {$\mathbf{w}_i$}  |- (aircraft.145);
	%
	% close outer velocity loop
	\node [block, below of=Kstar]  (Kvt) {$\frac{1}{s}K_v + K_{x_v}$};
	\draw [->] (aircraft.-10) -- ++(1,0) |- node[pos=0.75] {$\mathbf{v}_i$} (Kvt);
	\draw [->] (Kvt) -| node[pos=0.99] {$-$} (sumMiddleLoop);
	%|- ++(-2,0) node[pos=0.99] {$\tilde{\mathbf{v}}_i$} -|  node[pos=0.99] {$-$} (sumMiddleLoop.south);
	
	% close outer position loop
	\draw [->] (aircraft.10) -| ++(1.25,-3.5) -| node[pos=0.02] {$\mathbf{p}_i$} node[pos=0.99] {$-$} (sumOuter.south);

\end{tikzpicture}
\caption{Diagram on the controller whose gains will be tuned using structured H-infinity synthesis}
\label{fig:controlDiagram}
\end{figure}
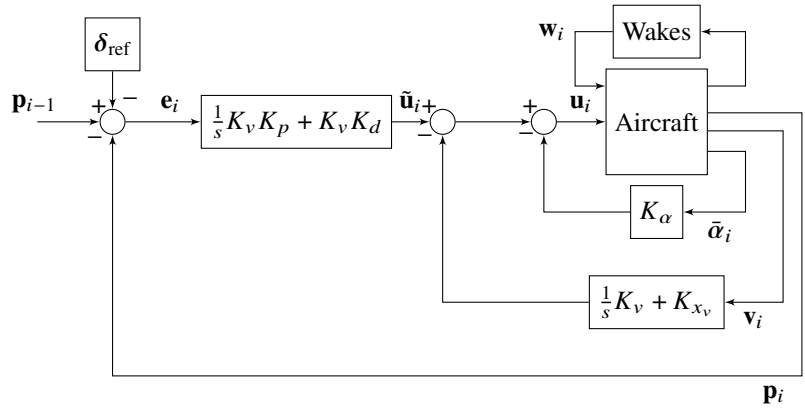

Fig. \ref{fig:controlDiagram} shows a diagram of the proposed control architecture, separating the different groups of control gains to emphasize their distinct roles in the control design.
The control surfaces and thrust are driven by the input $\control_i$, which is composed of proportional feedback gains $K_\alpha$ on the rotational states, plus proportional-integral (PI) controllers on both the separation error $\e_i$ and the deviation from nominal velocity $\v_i = \dot \p_i - \v_0$.
The rotational state feedback term $K_\alpha$ is responsible for stabilizing the attitude of the aircraft, while the PI term on the separation error should ensure that each aircraft accurately tracks the reference separation between the preceding aircraft.
\revjr{Finally, the feedback term on the deviation from trimmed velocity serves two purposes. It makes it possible to achieve string stability in the same way as a time headway, by eliminating one of the poles at the origin, and it regulates the steady-state velocity to ensure that all airplanes maintain the same velocity even when subjected to aerodynamic disturbances.}

The next step is to express the control objectives in terms of the proposed control gain matrices.
Recall that our linearized aircraft dynamics are given by:
\begin{equation*}
    \dot \x_i = A \x_i + B \u_i + \w_i,
\end{equation*}
where $\x$ denotes the 12-dimensional state, $\u$ the four-dimensional input, and $\w$ a 12-dimensional disturbance modeling the wake effects.
In the notation above, we have decomposed the state as $\x_i := [\p_i^\top\,\, \v_i^\top\,\, \bar{\a}_i^\top]^\top$, where $\bar{\a}_i := [\a_i^\top\,\,\dot\a_i^\top]^\top$ denotes the Euler angles and rotation rates grouped together (see Section \ref{sec:aircraft}).
The proposed controller can be written as
\begin{align*}
    \u_i = \left(\frac{1}{s}K_v K_p + K_v K_d\right)(\p_{i-1} - \p_i) - \left(\frac{1}{s}K_v + K_{x_v}\right)\tilde\v_i - K_\alpha \bar\a_i.
\end{align*}

Since $\u_i = \tilde{\u}_i - \left(\frac{1}{s}K_v + K_{x_v}\right)\tilde\v_i - K_\alpha \bar\a_i$, we have
\begin{align*}
    \tilde{\u}_i &=  \left(\frac{1}{s}K_v K_p + K_v K_d\right)(\p_{i-1} - \p_i) \\
    &= C(s) (\p_{i-1} - \p_i),
\end{align*}
where the outer-loop controller is given by $C(s) =  \left(\frac{1}{s}K_v K_p + K_v K_d\right)$.

Let $P(s)$ denote the open-loop transfer function for the aircraft from $\u_i$ to $\p_i$.
Then let $\bar{P}(s) = C_p(I + (\frac{1}{s}K_v + K_{x_v})C_v + K_\alpha C_\alpha)^{-1}P$ denote the transfer function from $\tilde{\u}_i$ to $\p_i$, where $C_p$, $C_v$, and $C_\alpha$ are matrices that isolate the position, velocity, and angular states, respectively, from the full state vector $\x_i$ (e.g. $\p_i = C_p \x_i$).
To analyze string stability, we are interested in $T(s)$ (the complementary sensitivity transfer matrix) from $\p_{i-1}$ to $\p_i$:
\begin{align*}
   \p_i &= \bar P(s) \tilde{\u}_i \\
   \p_i &= \bar P(s) C(s) (\p_{i-1} - \p_i) \\
   \p_i &= (I + \bar P(s) C(s))^{-1} \bar P(s) C(s) \p_{i-1} \\
   T(s) &=  (I + \bar P(s) C(s))^{-1} \bar P(s) C(s).
\end{align*}
Recall that to achieve string stability in the formation, $T(s)$ must satisfy $\sup_\omega \bar{\sigma}[T(j\omega)] \leq 1$.
We can equivalently express this as the H-infinity constraint $||T(s)||_\infty \leq 1$.
While there exist multiple methods for H-infinity synthesis, since we already have a target control architecture, we choose the structured H-infinity design method of \cite{gahinet2011hinf}. This approach involves using non-smooth optimization on a set of tunable system parameters after transforming the system into a canonical feedback form.

\begin{figure}[ht]
    \begin{center}
    \includegraphics[width=0.7\linewidth]{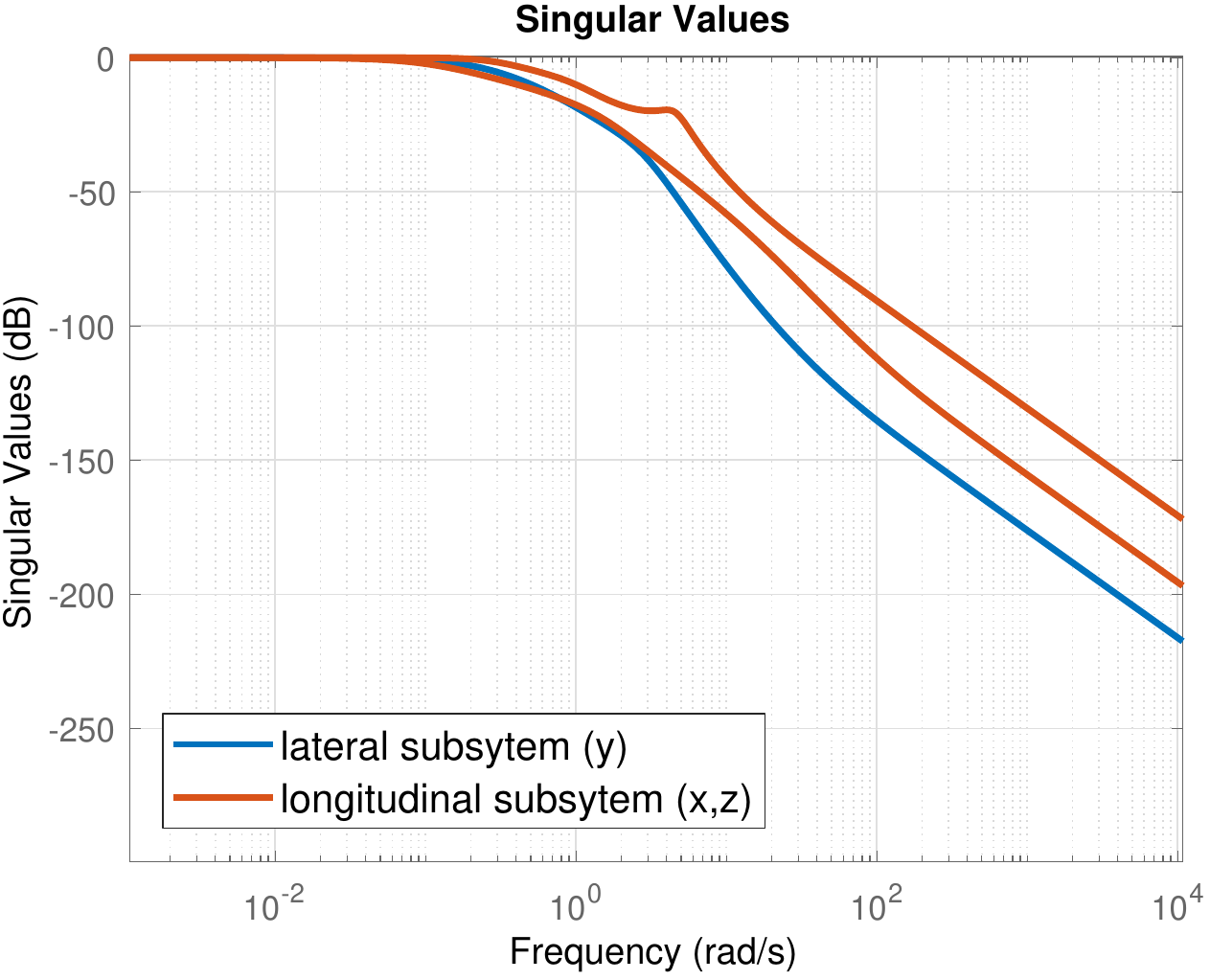}
    \caption{Singular values of the complementary sensitivity transfer matrix $T(j\omega)$ using the structured H-infinity control design} \label{fig:sv}
    \end{center}
\end{figure}

We then solve the H-infinity optimization to ensure that $||H||_\infty \leq 1$ using the MATLAB function \code{hinfstruct} from the Robust Control Toolbox.
We choose the diagonal elements of the control gain matrices $K_p$ and $K_d$ as tunable parameters and initializing them to arbitrary values, ensuring that the system remains stable.
In addition to the string stability constraint, we provide two additional performance specifications that restrict the search space of the optimization algorithm: (i) we set the minimum decay rate to .08 $\text{s}^{-1}$ to ensure a sufficiently fast response time, and (ii) we set the maximum frequency to 50 rad/s to prevent unreasonably large control gains.
Figure \ref{fig:sv} shows the singular values of $T(s)$ for the resulting controller, where we see that the system is indeed string stable.

\begin{figure}[ht]
    \begin{center}
    \includegraphics[width=0.7\linewidth]{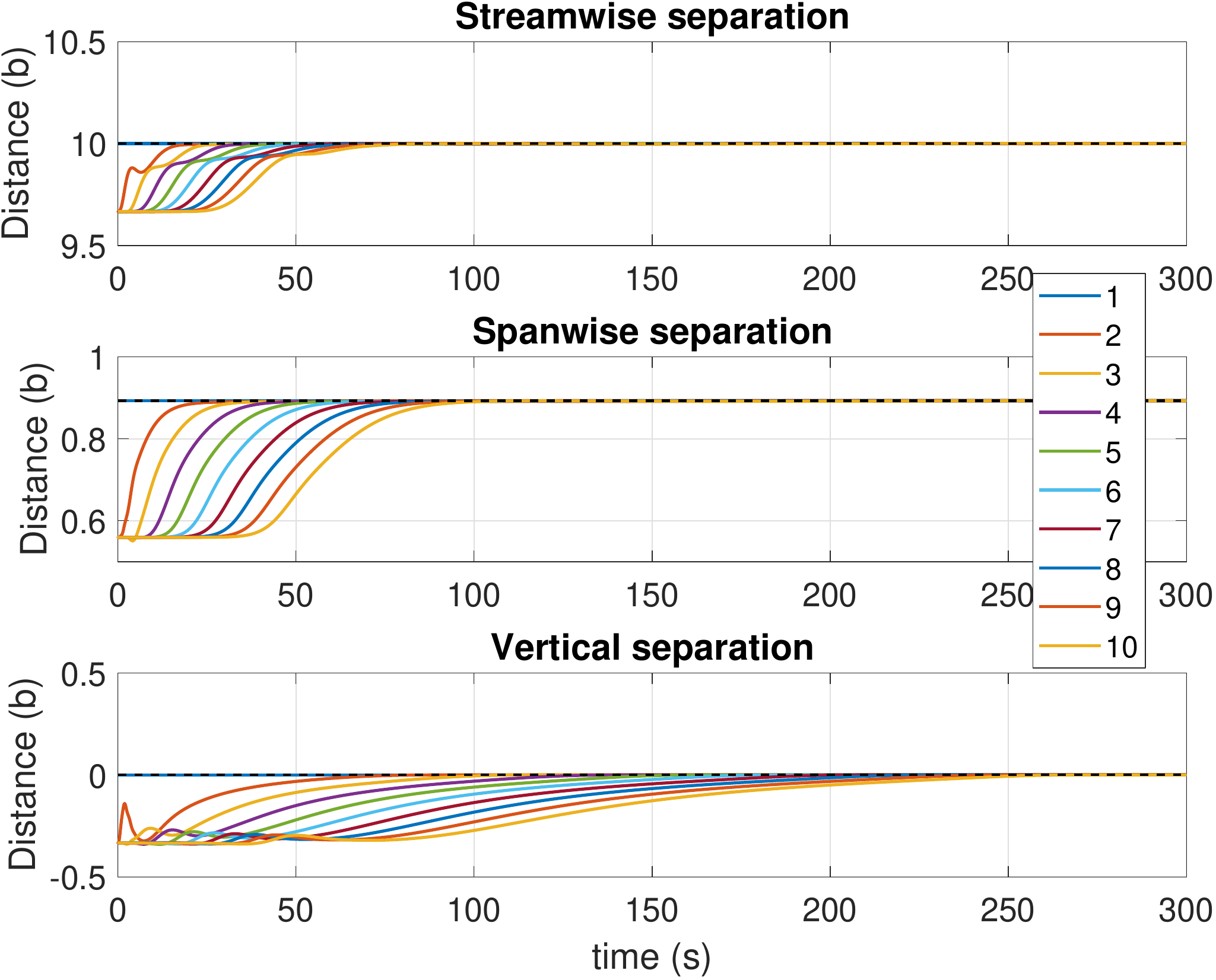}
    \caption{Using the proposed controller, each aircraft in the formation converges to the target offset vector while exhibiting string stable behavior.} \label{fig:pdlqrsim}
    \end{center}
\end{figure}

\begin{figure}[ht]
    \centering
    \includegraphics[width=0.7\linewidth]{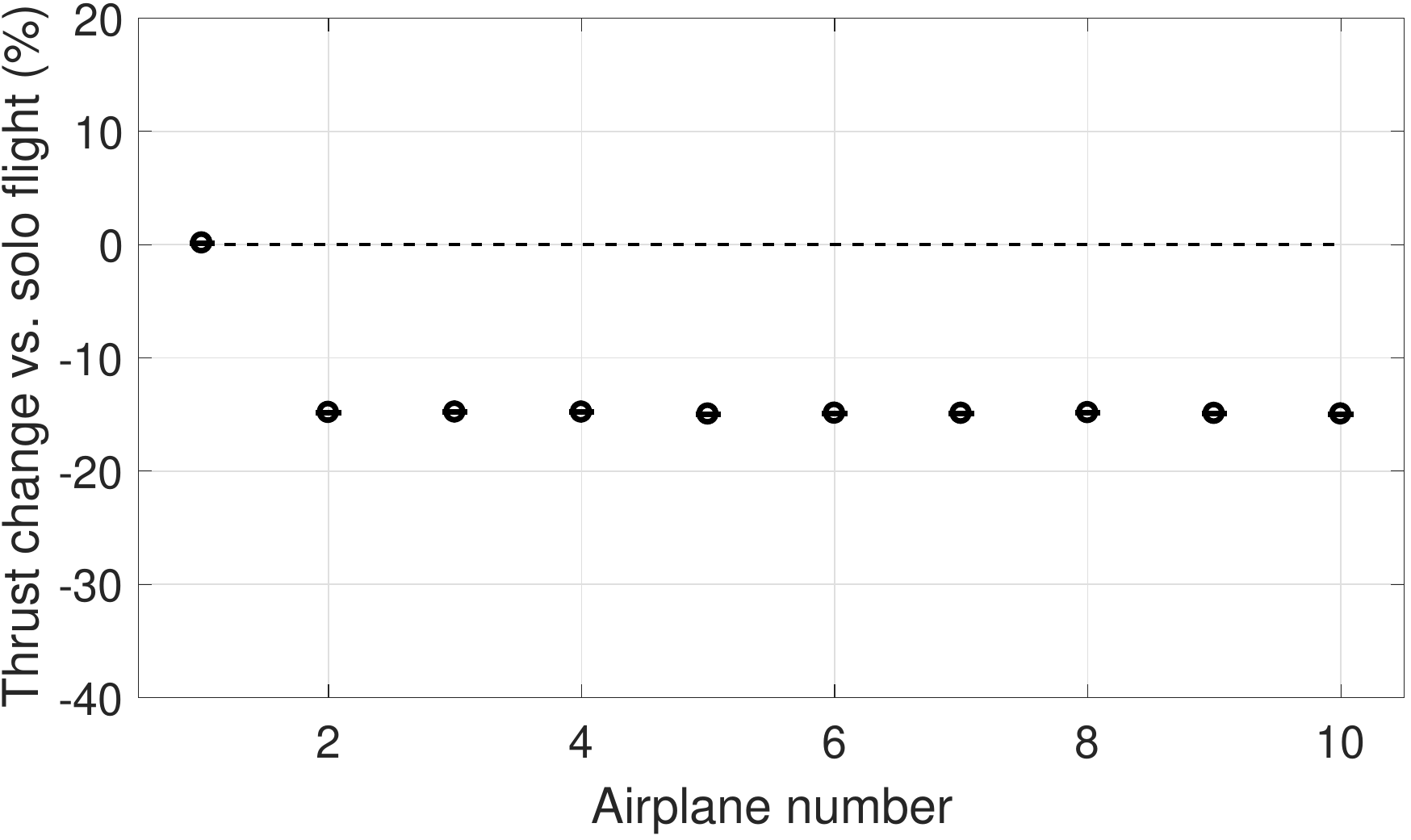}
    \caption{Thrust change compared to solo flight in 2\% turbulence intensity.}
	\label{fig:energy_pdlqr}
\end{figure}

The response of this controller to an initial perturbation is shown in Fig. \ref{fig:pdlqrsim}, demonstrating both string stability and convergence to the target positions.
The behavior of the formation and corresponding energy savings in 2\% turbulence intensity are shown in Fig. \ref{fig:pdlqrsim} and \ref{fig:energy_pdlqr}.
The steady state thrust is averaged over the final 30 seconds of flight, and the corresponding standard deviations are indicated by the error bars. Although the system is perturbed by the turbulence, these disturbances are no longer amplified along the formation, so there is minimal degradation of the drag reduction.

\section{Conclusions}

Groups of airplanes can save significant amounts of energy by flying in formations that take advantage of the beneficial upwash regions in the wakes of preceding airplanes.
\rev{To achieve good performance, each following airplane must accurately position itself in the wake of the aircraft immediately in front.
At the same time, disturbances induced by wind gusts and wake effects must not be amplified along the formation, i.e. the formation must be string stable.
In this paper, we have demonstrated that achieving both string stability and good tracking performance in the presence of aerodynamic disturbances is a challenge for aircraft formations, although for slightly different reasons than in other domains such as automobile platoons.
Nevertheless, we have shown one method for designing a controller using local velocity feedback that achieves both objectives.
Simulations show that such a controller can allow for arbitrarily long aircraft formations with energy savings of approximately 15\% for each following aircraft.}

We hope to elaborate on this work in several ways. 
\rev{Adding actuator dynamics to the aircraft models could reveal whether this impacts the string stability properties of the formation. 
The flow disturbance models can also be improved, for example by deforming the wakes as they move downstream. 
Additionally, the turbulence model can be enhanced by including variations in the $y$ and $z$ directions. 
Finally, we hope to extend this work to the case of tracking the actual wake position rather than the position of the immediate leader.
Including wake estimation, as described in \cite{caprace2019wake}, will lead to several new challenges, but has the potential to further improve the performance.}

\clearpage

\section*{Appendix}
\setcounter{subsection}{0} % otherwise it starts at B, for some reason
\subsection{Aircraft dimensions}
\begin{table}[ht]
    \centering
    \begin{tabular}{c c}
        \hline 
        Property & Value \\
        \hline 
        Mass & 80,000 kg \\ 
        Wingspan & 34.1 m\\ 
        Mean chord & 3.6 m\\ 
        Cruise speed & 230 m/s\\ 
        Air density & 0.458 kg/m$^3$\\ 
        Tail span & 12.5 m\\ 
        Vertical tail span & 6.2 m\\ 
        Trimmed thrust & $5.02\times10^{4}$ N\\ 
        Zero-lift drag coefficient $C_{D,0}$ & 0.03\\ 
        Wake circulation, $\Gamma$ & 278 m$^2$/s\\
        \hline
    \end{tabular}
    \caption{Aircraft dimensions and cruise flight properties, from Colognesi\cite{colognesi}}
    \label{tab:aircraftDimensions}
\end{table}

\subsection{Potential flow velocities}
\subsubsection{Vortex filament}
The velocity at point $\mathbf{p}$ due to a straight finite vortex filament from $\mathbf{p}_1$ to $\mathbf{p}_2$ with circulation $\Gamma$ is:
\begin{subequations}
\begin{align}
	\mathbf{r}_0 &=\mathbf{p}_1-\mathbf{p}_2 \\
	\mathbf{r}_1 &=\mathbf{p}-\mathbf{p}_2 \\
	\mathbf{r}_2 &=\mathbf{p}-\mathbf{p}_1  \\
	U_{vf} &=   \frac{\Gamma}{4 \pi } \left( \frac{ \mathbf{r}_1 \times \mathbf{r}_2 }{r_c^2+|\mathbf{r}_1 \times \mathbf{r}_2|^2 } \right)  \left( \frac{\mathbf{r}_0 \cdot \mathbf{r}_1}{|\mathbf{r}_1|^2} -\frac{\mathbf{r}_0 \cdot \mathbf{r}_2}{|\mathbf{r}_2|^2}  \right). \label{eq:vortexFilament}
\end{align}
\end{subequations}
For the VLM computations, $r_c=$ 0.1mm to remove the singularity.

\subsubsection{Horseshoe vortex}
The velocity due to an aircraft and its wake has velocity $U_{HH}$ from the bound vortex ``head'' and $U_{HL}$ and $U_{HR}$ from the semi-infinite left and right wake vortex ``legs''. The two legs are separated by a horizontal distance of $d=b \pi/4$. For this notation, the head of the vortex is centered at $\mathbf{p}_v = [x_v,y_v,z_v]$, the left and right corners are at $\mathbf{p}_L=\mathbf{p}_v- \mathbf{\hat{j}}d/2$ and $\mathbf{p}_R =\mathbf{p}_v + \mathbf{\hat{j}}d/2$, and the velocity is computed at $\mathbf{p} = [x,y,z]$. The core size of the vortex is $r_c=0.05b$. The velocity field due to such a vortex is:
\begin{subequations}
\begin{align}
	\mathbf{r}_0 &=\mathbf{p}_R-\mathbf{p}_L \\
	\mathbf{r}_1 &=\mathbf{p}-\mathbf{p}_L \\
	\mathbf{r}_2 &=\mathbf{p}-\mathbf{p}_R  \\
	U_{HH} &=   \frac{\Gamma}{4 \pi } \left( \frac{ \mathbf{r}_1 \times \mathbf{r}_2 }{r_c^2+|\mathbf{r}_1 \times \mathbf{r}_2|^2 } \right)  \left( \frac{\mathbf{r}_0 \cdot \mathbf{r}_1}{|\mathbf{r}_1|^2} -\frac{\mathbf{r}_0 \cdot \mathbf{r}_2}{|\mathbf{r}_2|^2}  \right) \\
	U_{HL} &=   \frac{\Gamma}{4 \pi } \left(  \frac{-(z-z_L)\mathbf{\hat{j}}+(y-y_L)\mathbf{\hat{k}}}{r_c^2+(y-y_L)^2+(z-z_L)^2  } \right) \left( 1-\frac{x-x_L}{ \sqrt{ (x-x_L)^2+(y-y_L)^2+(z-z_L)^2 } }  \right) \\
	U_{HR} &=   \frac{\Gamma}{4 \pi } \left(  \frac{(z-z_R)\mathbf{\hat{j}}-(y-y_R)\mathbf{\hat{k}}}{r_c^2+(y-y_R)^2+(z-z_R)^2  } \right) \left( 1-\frac{x-x_R}{ \sqrt{ (x-x_R)^2+(y-y_R)^2+(z-z_R)^2 } }  \right) \\
	U_H &= U_{HH}+U_{HL}+U_{HR} \label{eq:wakeVel}
\end{align}
\end{subequations}

\subsection{Linearized dynamics}
\label{app:systemDynamics}

The nomenclature for the state and control variables is displayed in Tables \ref{tab:variableNamesState} and \ref{tab:variableNamesControl}. The state and control vectors are presented in Equations \ref{eq:stateControlVectors}, and the linearized aircraft dynamics in Equation \ref{eq:systemMatrices}.

\begin{table}[ht]
    \centering
    \begin{tabular}{c|c | c}
        \hline
        Symbol & Meaning & Units \\
        \hline
        $x$ & Streamwise position & m\\
        $y$ & Lateral position & m\\
        $z$ & Vertical position, positive down & m\\
        $\dot{x}$ & Streamwise velocity & m/s\\
        $\dot{y}$ & Spanwise velocity & m/s\\
        $\dot{z}$ & Vertical velocity & m/s\\
        $\phi$ & Roll & radians\\
        $\theta$ & Pitch & radians\\
        $\psi$ & Yaw & radians\\
        $\dot{\phi}$ & Roll rate & radians/s\\
        $\dot{\theta}$ & Pitch rate & radians/s\\
        $\dot{\psi}$ & Yaw rate & radians/s\\
        \hline
    \end{tabular}
    \caption{Variables in the state vector}
    \label{tab:variableNamesState}
\end{table}

\begin{table}[ht]
    \centering
    \begin{tabular}{c|c | c}
        \hline
        Symbol & Meaning & Units \\
        \hline
        $\Delta T$ & Thrust (change from cruise)& N\\
        $\Delta_a$ & Aileron deflection & radians\\
        $\Delta_e$ & Elevator deflection & radians\\
        $\Delta_r$ & Rudder deflection & radians\\
        \hline
    \end{tabular}
    \caption{Variables in the control vector}
    \label{tab:variableNamesControl}
\end{table}

The state and control vectors, separated into longitudinal and lateral modes, are:
\setcounter{MaxMatrixCols}{20}
\begin{subequations}
\begin{align}
    \mathbf{x} &= 
    \begin{pmatrix}
        x & y & z & \dot{x} & \dot{y} & \dot{z} & \phi & \theta & \psi & \dot{\phi} & \dot{\theta} & \dot{\psi} \\
    \end{pmatrix}^T, \\
    \mathbf{x}_{long} &= 
    \begin{pmatrix}
        x &  z & \dot{x} &  \dot{z} &  \theta & \dot{\theta} \\
    \end{pmatrix}^T, \\
    \mathbf{x}_{lat} &= 
    \begin{pmatrix}
        y & \dot{y} & \phi  & \psi & \dot{\phi}  & \dot{\psi} \\
    \end{pmatrix}^T,\\
    \mathbf{u} &= 
    \begin{pmatrix}
        T & \Delta_a & \Delta_e  & \Delta_r \\
    \end{pmatrix}^T.
\end{align}
\label{eq:stateControlVectors}
\end{subequations}
The linearized aircraft dynamics, similarly separated, are:
\begin{subequations}
\begin{align}
    A_{long} &=
    \begin{pmatrix}
        0 & 0 & 1 & 0 & 0 & 0 \\
        0 & 0 & 0 & 1 & 0 & 0 \\
        0 & 0 & -5.45\text{e}{-3} & 3.61\text{e}{-2} & -1.51 & -6.42\text{e}{-2} \\
        0 & 0 & -8.52\text{e}{-2} & -0.445 & -102 & 227 \\
        0 & 0 & 0 & 0 & 0 & 1 \\
        0 & 0 & 0 & -4.18\text{e}{-2} & -9.62 & -0.960 \\
    \end{pmatrix}, \\
    A_{lat} &=
    \begin{pmatrix}
        0 & 1 & 0 & 0 & 0 & 0 \\
        0 & -3.57\text{e}{-2} & 9.81 & 8.22 & -0.167 & -230 \\
        0 & 0 & 0 & 0 & 1 & 0 \\
        0 & 0 & 0 & 0 & 0 & 1 \\
        0 & -1.10\text{e}{-2} & 0 & 2.52 & -0.395 & 0.193 \\
        0 & 6.29\text{e}{-3} & 0 & -1.45 & -4.76\text{e}{-3} & -0.135
    \end{pmatrix}, \\
%\end{align}
%\end{subequations}
%\begin{subequations}
%\begin{align}
    B_{long} &=
    \begin{pmatrix}
        0 & 0 & 0 & 0 \\
        0 & 0 & 0 & 0 \\
        1.25\text{e}{-5} & 0 & -0.138 & 0 \\
        0 & 0 & -7.20 & 0 \\
        0 & 0 & 0 & 0 \\
        0 & 0 & -3.50 & 0   \\
    \end{pmatrix}, \\
    B_{lat} &=
    \begin{pmatrix}
        0 & 0 & 0 & 0 \\
        0 & 0.487 & 0 & 4.59 \\
        0 & 0 & 0 & 0 \\
        0 & 0 & 0 & 0 \\
        0 & 1.08 & 0 & 0.418 \\
        0 & -1.82\text{e}{-2} & 0 & -0.960 \\
    \end{pmatrix} .
\end{align}
\label{eq:systemMatrices}
\end{subequations}
% done, with no noticeable difference \todo{remove very small values from A,B?}

% LQR gains (K_x)

\begin{landscape}
\subsection{Controller specifications} \label{app:controlSpecs}

\subsubsection{LQR}

This section displays the control matrices that were found by optimizing the LQR weight matrices for disturbance attenuation from leading to following aircraft. \done{PC: add some text to introduce the LQR matrices}
\revThree{The controller is given by $\u = -K_x \x$ where}
\newcommand{\timesb}{\hspace{-0.15em} \times \hspace{-0.15em}}
\begin{subequations}
\footnotesize
\begin{align}
    K_x &=
    \begin{pmatrix}
        2.23\text{e}{+4} & -3.48\text{e}{-8} & -916 & 5.93\text{e}{+4} & 1.05\text{e}{-8} & -177 & 8.25\text{e}{-7} & 5.54\text{e}{+4} & 3.98\text{e}{-6} & 3.54\text{e}{-7} & 1.19\text{e}{+4} & 3.05\text{e}{-7} \\
        0 & 7.75\text{e}{-3} & 0 & 0 & 4.25\text{e}{-2} & -3.91\text{e}{-10} & 0.751 & 9.24\text{e}{-8} & 6.65  & 0.828 & 3.17\text{e}{-9} & -0.740 \\ 
        9.16\text{e}{-4} & 0 & 4.45\text{e}{-3} & -7.74\text{e}{-4} & 0 & 1.98\text{e}{-2} & 0 & -4.70  & 0 & 0 & -0.167 & 0 \\
        0 & 9.70\text{e}{-3} & -3.45\text{e}{-10} & 0 & 6.63\text{e}{-2} & -1.07\text{e}{-9} & 0.192 & 2.52\text{e}{-7} & 1.10  & 2.52\text{e}{-3} & 7.24\text{e}{-9} & -4.96  
    \end{pmatrix}, 
\end{align}
\end{subequations}

% LQR plus integral gains 
\subsubsection{LQR plus integral}

\revThree{For the LQR plus integral controller, the state is augmented with the integral of the position error in the three spatial dimensions, that is, $\bar \x := [\x^\top \,\, \x_{int}^\top]^\top$ and $\dot \x_{int} = \e$.
The controller is then given by $\u = -K_{\bar x} \bar\x$, where}
\begin{subequations}
\footnotesize
\begin{align}
    K_{\bar x} &=
    \begin{pmatrix}
        3.04\text{e}{+4} & -6.24\text{e}{-7} & -3.27\text{e}{+3} & 6.97\text{e}{+4} & 3.27\text{e}{-8} & -3.83\text{e}{+3} & 2.02\text{e}{-6} & 9.07\text{e}{+5} & 8.93\text{e}{-6} & 1.49\text{e}{-7} & 2.59\text{e}{+4} & 9.15\text{e}{-7} & 3.14\text{e}{+3} & -1.23\text{e}{-7} & -413 \\
        0 & 3.46\text{e}{-2} & 5.65\text{e}{-10} & 0 & 5.07\text{e}{-2} & 1.08\text{e}{-9} & 0.770 & -2.53\text{e}{-7} & 6.77  & 0.834 & -4.71\text{e}{-9} & -1.04  & 0 & 1.05\text{e}{-2} & 0 \\
        2.50\text{e}{-3} & 0 & 1.65\text{e}{-2} & -1.16\text{e}{-3} & 0 & 4.44\text{e}{-2} & 0 & -10.4 & 0 & 0 & -0.283 & 0 & 1.85\text{e}{-4} & 0 & 1.40\text{e}{-3} \\
        0 & 7.71\text{e}{-2} & -2.71\text{e}{-10} & 0 & 8.63\text{e}{-2} & -5.19\text{e}{-10} & 0.231 & 1.21\text{e}{-7} & 1.32  & 1.13\text{e}{-2} & 2.28\text{e}{-9} & -5.70  & 0 & 3.14\text{e}{-2} & 0 \\
    \end{pmatrix}.
\end{align}
\end{subequations}
\end{landscape}

\subsubsection{Proposed controller gains}

\revThree{The proposed controller designed by structured H-infinity synthesis is given by}
\begin{align*}
    \u_i = \left(\frac{1}{s}K_v K_p + K_v K_d\right)(\p_{i-1} - \p_i) - \left(\frac{1}{s}K_v + K_{x_v}\right)\tilde\v_i - K_\alpha \bar\a_i,
\end{align*}
where
\begin{subequations}
\footnotesize
\begin{align}
    K_\alpha &=
   \left(\begin{array}{ccccccc} -2.302\,{10}^5 & 9.372\,{10}^{-5} & 5.411\,{10}^7 & 0.0007509 & 4.229\,{10}^{-5} & 9.616\,{10}^5 & -0.0007284\\ -9.863\,{10}^{-9} & 0.5396 & 2.309\,{10}^{-6} & 4.472 & 0.6881 & 3.655\,{10}^{-8} & 0.2467\\ 0.09398 & 0 & -22.63 & -3.302\,{10}^{-10} & 0 & -0.7453 & 2.694\,{10}^{-10}\\ -5.394\,{10}^{-8} & 0.1307 & 1.262\,{10}^{-5} & 1.076 & -0.01545 & 1.878\,{10}^{-7} & -3.75 \end{array}\right)
\end{align}
\end{subequations}

\begin{subequations}
\begin{align}
    K_v &=\left(\begin{array}{ccc} 84677.0 & -6.893\,{10}^{-5} & -1.239\,{10}^5\\ -6.159\,{10}^{-10} & 0.009398 & -5.512\,{10}^{-9}\\ 0.005323 & 0 & 0.0291\\ -3.348\,{10}^{-9} & 0.03092 & -3.105\,{10}^{-8} \end{array}\right) \\
    K_p &= \left(\begin{array}{ccc} 0.2421 & 0 & 0\\ 0 & 0.1559 & 0\\ 0 & 0 & 0.07919 \end{array}\right) \\
    K_d &=\left(\begin{array}{ccc} 0.1006 & 0 & 0\\ 0 & 0.01063 & 0\\ 0 & 0 & 0.1746 \end{array}\right) \\
    K_{x_v} &=\left(\begin{array}{ccc} 1.318\,{10}^5 & 1.606\,{10}^{-5} & -2.302\,{10}^5\\ 1.067\,{10}^{-10} & 0.01954 & -9.863\,{10}^{-9}\\ -0.001378 & 0 & 0.09398\\ 5.834\,{10}^{-10} & 0.03872 & -5.394\,{10}^{-8} \end{array}\right)
\end{align}
\end{subequations}

\subsection{Transfer functions} \label{app:transferFunctions}

Open-loop lateral aircraft dynamics: from ailerons and rudder to lateral position:
\begin{align*}
    P_{a,y}(s) &= \frac{0.4868 s^4 + 4.247 s^3 + 13.97 s^2 + 1.324 s + 14.87}{s^6 + 0.5657 s^5 + 2.962 s^4 + 1.275 s^3 + 0.002584 s^2 + \num{5.131e-10} s} \\ ~&~ \\
    P_{r,y}(s) &= \frac{4.588 s^4 + 222.7 s^3 + 90.62 s^2 - 1.417 s - 17.81}{s^6 + 0.5657 s^5 + 2.962 s^4 + 1.275 s^3 + 0.002584 s^2 + \num{5.131e-10} s}.
\end{align*}

Open-loop transfer function composed of lateral aircraft dynamics and lateral LQR plus integral controller:
\begin{align*}
 P_{y,y}(s) &= \frac{0.1494\,s^4+7.048\,s^3+2.996\,s^2-0.03066\,s-0.4041}{s^7+0.5657\,s^6+2.962\,s^5+1.275\,s^4+0.002584\,s^3+5.131\,{10}^{-10}\,s^2}
\end{align*}

% \begin{align*}
%      P_{y,y}(s) &= \frac{1.405 s^6 + 1.963 s^5 + 5.244 s^4 + 4.655 s^3 + 2.625 s^2 + 0.6138 s + 0.06606}{s^7 + 0.5657 s^6 + 2.962 s^5 + 1.275 s^4 + 0.002622 s^3 + \num{1.859e-06} s^2} \\ ~&~ \\
%      P_{r,y}(s) &= \frac{6.134 s^6 + 11.04 s^5 + 8.135 s^4 - 0.1733 s^3 - 2.661 s^2 - 0.7347 s - 0.07916}{s^7 + 0.5657 s^6 + 2.962 s^5 + 1.275 s^4 + 0.002622 s^3 + \num{1.859e-06} s^2}.
% \end{align*}

% \frac{0.1494\,s^4+7.048\,s^3+2.996\,s^2-0.03066\,s-0.4041}{s^7+0.5657\,s^6+2.962\,s^5+1.275\,s^4+0.002584\,s^3+5.131\,{10}^{-10}\,s^2}

Complementary sensitivity function for LQR control applied to lateral aircraft dynamics:
\begin{align*}
    T_y(s) = \frac{0.02055 s^4 + 0.6886 s^3 + 1.203 s^2 + 0.8138 s + 0.4735}{s^6 + 4.199 s^5 + 11.61 s^4 + 14.65 s^3 + 10.13 s^2 + 4.519 s + 0.4735}.
\end{align*}

\section*{Funding Sources}
This project has received funding from the European Research Council (ERC) under the European Union’s Horizon 2020 research and innovation program (grant agreement No 725627) and from the French community of Belgium in the Joint Research Activity RevealFlight (convention no. 17/22-080).

\end{document}